\documentclass[twocolumn]{svjour3} 
\usepackage{graphicx}
\usepackage{latexsym}
\usepackage{amssymb}
\usepackage{subfigure}
\usepackage{fancyhdr}
\usepackage[usenames,dvipsnames]{pstricks}
\usepackage{epsfig}
\usepackage[english]{babel}
\usepackage{indentfirst}
\usepackage{subfigure}
\usepackage{xspace} 
\usepackage{amsmath,amssymb}

\makeatletter
\newcommand{\sech}{\mathop{\operator@font sech}}
\newcommand{\sign}{\mathop{\operator@font sign}}
\makeatother

\numberwithin{equation}{section}

\begin{document}
\title{Cross-diffusion systems for image processing: I. The linear case}
\author{A. Ara\'ujo         \and
        S. Barbeiro \and
E. Cuesta\and
A. Dur\'an
}
\authorrunning{Ara\'ujo {\em et al.}}
\institute{A. Ara\'ujo \at
              CMUC, Department of Mathematics, University of Coimbra, Portugal\\
              \email{alma@mat.uc.pt}           
           \and
           S. Barbeiro \at
              CMUC, Department of Mathematics, University of Coimbra, Portugal\\
\email{silvia@mat.uc.pt}
\and
E. Cuesta \at
              Department of Applied Mathematics,  University of
Valladolid,  Spain\\
\email{eduardo@mat.uva.es}
\and
A. Dur\'an \at
              Department of Applied Mathematics,  University of
Valladolid,  Spain\\
\email{angel@mac.uva.es}
}
\date{Received: date / Accepted: date}

\maketitle
\begin{abstract}
The use of cross-diffusion problems as mathematical models of different image processes is investigated. Here the image is represented by two real-valued functions which evolve in a coupled way, generalizing the {approaches} based on real and complex diffusion. The present paper is concerned with linear filtering. First, based on principles of scale invariance, a scale-space axiomatic is built. Then some properties of linear complex diffusion are generalized, with particular emphasis on the use of one of the components of the cross-diffusion problem for edge detection. The performance of the cross-diffusion approach is analyzed by numerical means in some one- and two-dimensional examples.
\keywords{Cross-diffusion \and Complex diffusion\and Image {denoising}}
\end{abstract}
%


\section{Introduction}
This paper is devoted to the study of cross-diffusion systems in image processing.
Cross-diffusion models consist of evolutionary systems of diffusion type for at least two real-valued functions, where the evolution of each function is not independent of the others. Their use is widespread in areas like population dynamics, (see Galiano et al. \cite{GalianoGJ2001,GalianoGJ2003} and references therein). In the case of image processing, {two} previous related approaches must be mentioned. The first one is the so-called complex diffusion (Gilboa et al. \cite{GilboaSZ2004}).  Here the image is represented by a complex function and a filtering process is governed by a (nonlinear in general) partial differential equation (PDE) problem of diffusion type with a complex diffusion coefficient. {The main} properties of this approach, that will be relevant to the present paper, are {briefly described in} the following (see Gilboa et al. \cite{GilboaZS2001,GilboaSZ2002,GilboaSZ2004} for details). {On the one hand, the use of complex diffusion to model the filtering process assumes a distribution of the image features between the real and the imaginary parts of a complex function which evolve in a coupled way.} This mutual influence is governed by the complex diffusion matrix. (In particular, the initial condition of the corresponding PDE is always a complex function with the real part given by the actual image and zero imaginary part.) A second point to be emphasized is that this distribution may give to the imaginary part of the image a particular role as edge detector. In the linear case, when the phase angle of the (constant) complex diffusion coefficient is small, the imaginary part approximates a smoothed Laplacian of the original image, scaled by time. This property (called small theta approximation) can be used in nonlinear models as well when, in order to control the diffusion and the edge detection, the corresponding diffusion tensor is taken to be dependent only on the imaginary part of the image, instead of on the size of the {image} gradients as in many nonlinear {real diffusion} models, {with the consequent computational advantages.}

%

Complex diffusion can indeed be rewritten as a cross-diffusion model for the real and imaginary parts of the image function. This approach {was} considered by Lorenz et al. \cite{LorenzBZ}, to analyze the existence of global solution of a related cross-diffusion system. In addition to the theoretical advantage of having a global solution, the numerical examples also suggest that the filtering process with the cross-diffusion system preserves better the textures of the image when compared to Perona-Malik models (Perona \& Malik \cite{PeronaM}).

The previous references were the starting point of our research on cross-diffusion systems as mathematical models for image filtering. The general purpose of the present paper and the subsequent paper devoted to the nonlinear case, is to extend the complex diffusion approach to more general cross-diffusion systems, analyzing how to generalize known properties and deriving new ones. We point out the main contributions of the present paper.

\begin{itemize}
\item Our first purpose is to study the cross-diffusion approach ({regarding the complex diffusion as a particular case}) as a scale-space representation, see e.~g. Iijima \cite{Iijima1963}, Witkin \cite{Witkin1983}, {Lindeberg \cite{Lindeberg1994,Lindeberg2011}, Florack \cite{Florack1997}, Duits et al. \cite{DuitsFGH2004}}. When the features of the image are distributed in two components ${\bf u}=(u,v)^{T}$, we assume that the linear filtering formulation is described by a matrix convolution of the form
\begin{eqnarray}
(K\ast{\bf u})({\bf x})=\begin{pmatrix}(k_{11}\ast u)({\bf x})+(k_{12}\ast v)({\bf x})\\
(k_{21}\ast u)({\bf x})+(k_{22}\ast v)({\bf x})\end{pmatrix},\label{cdl1}
\end{eqnarray}
where
\begin{eqnarray}
K=\begin{pmatrix} k_{11}&k_{12}\\k_{21}&k_{22}\end{pmatrix},\label{cdl1b}
\end{eqnarray}
is the matrix kernel. For simplicity, and when no confusion is possible, $\ast$ will denote both the matrix-vector operation on the left hand side of (\ref{cdl1}) and the usual convolution of the entries on the right hand side,
\begin{eqnarray*}
(f\ast g)({\bf x})=\int_{\mathbb{R}^{2}}f({\bf x}-{\bf y})g({\bf y})d{\bf y}.\label{la1}
\end{eqnarray*}

In scale-space theory, two main formulations can be {considered}, see e.~g. Weickert et al. \cite{WeickertII,WeickertII2} and Lindeberg \cite{Lindeberg1997}. The first one is based on causality (Koenderink \cite{Koenderink1984}). This principle can be interpreted as non-enhancement of local extrema. The characterization of those linear filters in the continuous, one-dimensional case satisfying this property was given by Lindeberg in \cite{Lindeberg1990}. They consist of compositions of two types of  scale-space kernels: Gaussian kernels (some already observed by Witkin in \cite{Witkin1983}) and truncated exponential functions.

Our approach to build the scale-space axiomatics {makes use of a formulation} based on the principles of scale invariance (Iijima \cite{Iijima1963}, Florack et al. \cite{FloracktRHKV1992}, Pauwels et al. \cite{PauwelsVFN1995}, Weickert et al. \cite{WeickertII2}; see Lindeberg  \cite{Lindeberg1997} for a comparison and synthesis of both theories). This theory characterizes those linear filtering kernels satisfying recursivity, grey-level shift, rotational and scale invariance. Our first contribution in this paper generalizes that of the scalar case (that is, when the image is represented by only one real-valued function of two space variables) in the sense that the kernels $K=K_{t}$, being $t$ the scale parameter, have a Fourier matrix representation of the form
\begin{eqnarray}
\widehat{K}({\bf \xi},t)=e^{-t|{\bf \xi}|^{p}d},\quad p>0,\label{cdl2}
\end{eqnarray}
where $\widehat{K}({\bf \xi},t)$ stands for the $2\times 2$ matrix with entries $\widehat{k}_{ij}({\bf \xi},t),i,j=1,2$, $d=(d_{ij})_{i,j=1,2}$ is a $2\times 2$ positive definite matrix and $p$ is a positive constant. Additional properties which are analyzed here are the existence of infinitesimal generator (and consequently the alternative formulation of (\ref{cdl1}) as  an initial value problem of PDEs) and locality. The arguments in Pauwels et al. \cite{PauwelsVFN1995} and Duits et al \cite{DuitsFGH2004} will be here adapted to show that linear cross-diffusion filtering with convolution kernels of the form (\ref{cdl2})  admits an infinitesimal generator which is local if and only if $p$ is an even integer.
\item Since complex diffusion models can be written as cross-diffusion systems (for the real and imaginary parts of the image) we develop a generalization of {the most relevant} properties of the complex diffusion approach. They are related to the way how diffusion affects different features of the image, mainly the grey values and the edges. In this sense, the paper generalizes the small theta approximation with the aim of assigning to one of the components a similar role as edge detector to that of the imaginary part of the image in some cases of linear complex diffusion. The generalization is understood in the following terms. When the matrix $d$ approaches to {a suitable spectral} structure (see Section \ref{sec2}) then one of the components of the corresponding  filtered image by linear cross-diffusion behaves as the operator $A=-(-\Delta)^{p/2}$ applied to a smoothed version of the original image, determined by $A$ and the trace of $d$ (the sum of the diagonal entries of $d$). 

A second point of generalization concerns the initial distribution of the image in two components. In the complex diffusion case, linear (and nonlinear) filtering usually starts from the real noisy image written as complex function (that is, with zero imaginary part). This may be modified, {in the cross-diffusion approach}, by distributing the image in two components in different ways. This distribution affects the definition of relevant quantities for the problem in this context, such as the average grey level value.
\item The previous results are complemented with a computational study dealing with the performance of linear cross-diffusion filtering with kernels of the form (\ref{cdl2}). Our purpose here is to make a first approach to the behaviour of the models, by numerical means, that may be used in more exhaustive computational works in the future. The numerical experiments that we present are focused on the illustration of some properties of the models, as the generalized small theta approximation, and the influence of the initial distribution of the  {image features as well as the parameter $p$ and the matrix $d$ in (\ref{cdl2}). The numerical experiments suggest that when the eigenvalues of $d$ are real and different, the blurring effect in the filtering is delayed  when compared with} the linear complex case (where the eigenvalues are complex). The first choice leads to an improvement {in the quality of filtering, measured with the classical Signal-to-Noise-Ratio (SNR) and Peak-Signal-to-Noise-Ratio (PSNR) indexes, and edge detection} which is independent of the parameter $p$, that is, the local or nonlocal character of the infinitesimal generator.
\end{itemize}
The structure of the paper is as follows. Section \ref{sec2} is devoted to the theoretical analysis of linear cross-diffusion models. In Section \ref{sec21} the matrix convolution (\ref{cdl1}) is formulated as a scaling process and, in order to define the space where the convolution operators act, some hypotheses on the kernels are specified. The inclusion of a scaling enables to analyze, in Section \ref{sec22}, the formulation of some scale-space properties (grey level shift invariance, rotational invariance, semigroup property and scale invariance) under the cross-diffusion setting. The derivation of (\ref{cdl2}) is then given in Section \ref{sec23}, along with a discussion on the locality property. The generalization of the small theta approximation, in Section \ref{sec24}, completes the theoretical approach. Some of these properties are illustrated in the numerical study in Section \ref{sec3}, where the models are applied to one- and two-dimensional examples of filtering problems and their performance is analyzed in terms of the parameters of the kernels (\ref{cdl2}). Particular attention will be paid on the comparison with Gaussian smoothing and linear complex diffusion. The main conclusions and perspectives are outlined in Section \ref{sec4}. A technical result on the reduction of $2\times 2$ positive definite matrices, necessary for some proofs, is included as lemma in the Appendix.

The present paper will be followed by a {companion paper} focused on nonlinear cross-diffusion models (Ara\'ujo et al. \cite{ABCD2016_II}).

The following notation will be used throughout the paper. Grey-scale images will be represented in this linear case by real-valued mappings on $\mathbb{R}^{2}$. (The nonlinear case will be studied on bounded domains.) We denote by $H^{k}=H^{k}(\mathbb{R}^{2})$, where  $k\geq 0$ an integer, the Sobolev space of order $k$, with $H^{0}=L^{2}(\mathbb{R}^{2})$. The inner product in $H^{k}$ is denoted by $\langle\cdot,\cdot\rangle_{k}$ and the associated norm  by $||\cdot ||_{k}$. The space of integrable functions in $\mathbb{R}^{2}$ will be denoted as usual by $L^{1}=L^{1}(\mathbb{R}^{2})$, as well as the space of real-valued infinitely continuously differentiable functions in $\mathbb{R}^{2}$ by $C^{\infty}(\mathbb{R}^{2})$ and the space of continuous functions in $\mathbb{R}^{2}$ vanishing at infinity by $C_{0}(\mathbb{R}^{2})$. 
The vector and scalar real-valued functions on $\mathbb{R}^{2}$ will be distinguished by the use of boldface for the first. By $X$ we will denote the space where the scalar functions are, in such a way that a vector representation of the image will be defined in $X\times X$, with associated norm
\begin{eqnarray*}
||{\bf u}||_{X\times X}=\left(||u||_{X}^{2}+||v||_{X}^{2}\right)^{1/2},\quad {\bf u}=(u,v)^{T}.
\end{eqnarray*}
We will assume that $X=H^{0}$, although in some cases (which will be specified in the paper) $X=H^{k}$, for $k>0$, or $X=L^{1}\cap H^{0}$ will be considered. For $f\in X$, $\widehat{f}$ will denote the Fourier transform in $\mathbb{R}^{2}$,
\begin{eqnarray*}
\widehat{f}(\xi)=\int_{\mathbb{R}^{2}}f({\bf x})e^{-i{\bf x}\cdot \xi}d{\bf x}, \quad \xi\in\mathbb{R}^{2},
\end{eqnarray*}
where the dot $\cdot$ stands for the Euclidean inner product in $\mathbb{R}^{2}$ with $|\cdot |$ as the Euclidean norm. Finally, on the space of matrix kernels (\ref{cdl1b}) with $k_{ij}\in L^{1}, i,j=1,2$ we consider the norm
\begin{eqnarray}
\label{cdl1c}
||K||_{\ast}=\max_{i,j}\{||k_{ij}||_{L^{1}}\},
\end{eqnarray}
where $||\cdot ||_{L^{1}}$ denotes the $L^{1}-$norm.

\section{Linear cross-diffusion filtering}
\label{sec2}
\subsection{Formulation as scaling process}
\label{sec21}
In order to formalize (\ref{cdl1}) as a scale-space representation we introduce a family of convolution operators $\{\mathcal{K}_{t}:X\times X\rightarrow X\times X, t\geq 0\}$ in such a way that (\ref{cdl1}) is rewritten as
\begin{eqnarray}
{\bf u}({\bf x},t)=\mathcal{K}_{t}{\bf u}_{0}({\bf x})=(K(\cdot,t)*{\bf u}_{0})({\bf x}),\quad {\bf x}\in\mathbb{R}^{2},\label{cdl3}
\end{eqnarray}
where, from some original real-valued image $f\in X$, an initial vector field ${\bf u}_{0}({\bf x})=(u_{0}({\bf x}),v_{0}({\bf x}))^{T}\in X\times X$ is composed. Thus ${\bf u}({\bf x},t)=(u({\bf x},t),v({\bf x},t))^{T}$ stands for the grey-level value image at pixel ${\bf x}\in\mathbb{R}^{2}$ at the scale $t$. This is obtained from a convolution with a $2\times 2$ matrix kernel $K(\cdot,t)$ with entries $k_{ij}(\cdot,t), i,j=1,2$ such that the vector representation (\ref{cdl3}) is written as
\begin{eqnarray}
u({\bf x},t)&=&(k_{11}(\cdot,t)\ast u_{0})({\bf x})+(k_{12}(\cdot,t)\ast v_{0})({\bf x}),\nonumber\\
v({\bf x},t)&=&(k_{21}(\cdot,t)\ast u_{0})({\bf x})+(k_{22}(\cdot,t)\ast v_{0})({\bf x}),\nonumber\\
 &&t\geq 0, {\bf x}\in\mathbb{R}^{2}.\label{cdl4}
\end{eqnarray}
Concerning the kernels $k_{ij}, i.j=1,2$, a first group of hypotheses is assumed:
{
\begin{itemize}
\item[(H1)] $k_{ij}(\cdot,t)\in L^{1}, \widehat{k}_{ij}(\cdot,t)\in L^{1}, i,j=1,2, t> 0$.
\item[(H2)] For each ${\bf x}\in\mathbb{R}^{2}, i,j=1,2$, $k_{i,j}({\bf x},\cdot):(0,\infty)\rightarrow\mathbb{R}$ is continuous.
\end{itemize}
Note that since $\widehat{k}_{ij}(\cdot,t)\in L^{1}$ then $k_{ij}(\cdot,t): \mathbb{R}^{2}\rightarrow \mathbb{R}$ is continuous and bounded.
}
\begin{remark}
Hypotheses (H1), (H2) will be required for technical reasons (definition of convolution, inversion of Fourier transform) and also when some scale-space properties are imposed in (\ref{cdl3}). The satisfaction of these properties will require additional assumptions that will be specified in each case in Section \ref{sec22}.
\end{remark}
\begin{remark}
In a similar way to the scalar case ( Pauwels et al. \cite{PauwelsVFN1995}, Weickert et al. \cite{WeickertII2}, Lindeberg  \cite{Lindeberg1997}), it is not hard to see that the convolution kernel formulation (\ref{cdl3}) can be derived from the assumptions of linear integral operators  (in the matrix sense) {$\mathcal{K}_{t}, t\geq 0$ with matrix kernels $K_{t}$ such that}
\begin{eqnarray*}
\mathcal{K}_{t}{\bf f}({\bf x})=\int_{\mathbb{R}^{2}}K_{t}({\bf x},{\bf y}){{\bf f}}({\bf y})d{\bf y},\quad {\bf x}\in\mathbb{R}^{2},t\geq 0,
\end{eqnarray*}
and {satisfying the} translation invariance
\begin{eqnarray*}
K_{t}({\bf x}-{\bf y},\cdot)=K_{t}({\bf x},{\bf y}+\cdot),\quad {\bf x},{\bf y}\in\mathbb{R}^{2}, t\geq 0.
\end{eqnarray*}
\end{remark}
\begin{remark}
The linear complex diffusion with coefficient $c=re^{i\theta}$ can be written as a convolution (Gilboa et al. \cite{GilboaZS2001,GilboaSZ2002,GilboaSZ2004}) 
 \begin{eqnarray}
 I({\bf x},t)=(h(\cdot,t)\ast I_{0})({\bf x}),\quad {\bf x}\in\mathbb{R}^{2},\label{cdl5}
 \end{eqnarray}
 with kernel
 \begin{eqnarray}
 &&h({\bf x},t)=g_{\sigma}e^{i\alpha({\bf x},t)},
g_{\sigma}({\bf x})=\frac{1}{2\pi \sigma^{2}}e^{-\frac{|{\bf x}|^{2}}{2\sigma^{2}}},\nonumber\\
&&\sigma(t)=\sqrt{\frac{2tr}{\cos\theta}},\quad \alpha({\bf x},t)=\frac{|{\bf x}|^{2}\sin{\theta}}{4tr}.\label{cdl6}
 \end{eqnarray}
 If $I_{0}=I_{0R}+iI_{0I}, I=I_{R}+iI_{I}$ then (\ref{cdl5}) can be formulated as (\ref{cdl4}) for $u=I_{R}, u_{0}=I_{0R}, v=I_{I}, v_{0}=I_{0I}$ and $k_{11}({\bf x},t)=k_{22}({\bf x},t)=h_{R}({\bf x},t), k_{12}({\bf x},t)=-k_{21}({\bf x},t)=h_{I}({\bf x},t)$, where $h_{R}, h_{I}$ stand, respectively for the real and imaginary parts of $h$ in (\ref{cdl6}).
\end{remark}
\subsection{Scale-space properties}
\label{sec22}
As mentioned in the Introduction, the image representation as a scale-space will be here analyzed by using the principles of scale invariance. The purpose is then to characterize those matrix kernels $K(\cdot,t), t\geq 0$ in such a way that (\ref{cdl3}) satisfies shift invariance, rotational invariance, recursivity (semigroup property) and scale invariance. These four properties will be introduced in the context of cross-diffusion formulation and the requirements for the kernels to satisfy them will be imposed.
\subsubsection{Grey-level shift invariance}
As in the scalar case, we assume that 
\begin{itemize}
\item[(H3)]
The matrix kernel $K(\cdot,t), t>0$ is \lq mass-preserving\rq , Pauwels et al. \cite{PauwelsVFN1995}, that is 
\begin{eqnarray}
\widehat{k}_{ii}({\bf 0},t)&=&\int_{\mathbb{R}^{2}}k_{ii}({\bf x},t)d{\bf x}=1,\quad i=1,2,\nonumber\\
\widehat{k}_{ij}({\bf 0},t)&=&\int_{\mathbb{R}^{2}}k_{ij}({\bf x},t)d{\bf x}=0, \quad i\neq j.\label{cdl7}
\end{eqnarray}
\end{itemize}
Then for any constant signal $C\in\mathbb{R}^{2}$ and $t\geq 0$ we have $K(\cdot,t)\ast C=C$ and therefore grey-level shift invariance holds:
\begin{eqnarray*}
K(\cdot,t)\ast ({\bf f}+{\bf C})=K(\cdot,t)\ast{\bf f}+{\bf C},\quad {\bf C}\in\mathbb{R}^{2}, t\geq 0,
\end{eqnarray*}
for any input image ${\bf f}\in X\times X$. Hypothesis (H3) has an additional consequence:
\begin{lemma}
\label{lem1}
Assume that (H1)-(H3) hold.
For ${\bf u}=(u,v)\in L^{1}\times L^{1}$ we define 
${\bf M}({\bf u})=(m(u),m(v))^{T}$, where
\begin{eqnarray*}
m(f)=\left(\int_{\mathbb{R}^{2}}f({\bf x})d{\bf x}\right),\quad f\in L^{1}.
\end{eqnarray*}
If ${\bf u}_{0}\in L^{1}\times L^{1}$ and ${\bf u}(\cdot,t), t\geq 0$ satisfies (\ref{cdl3}) then
\begin{eqnarray}
{\bf M}({\bf u}(\cdot,t))=M({\bf u}_{0}),\quad t\geq 0.\label{cdl8}
\end{eqnarray}
\end{lemma}

{\em Proof}. Note that
for any $f\in L^{1}(\mathbb{R}^{2})$ and $i,j=1,2$,
\begin{eqnarray}
\int_{\mathbb{R}^{2}}\left(k_{ij}(\cdot,t)\ast f\right)({\bf x})d{\bf x}=\int_{\mathbb{R}^{2}}\int_{\mathbb{R}^{2}}k_{ij}({\bf x-y},t)f({\bf y})d{\bf y}d{\bf x}&&\nonumber\\
=
\int_{\mathbb{R}^{2}}\left(\int_{\mathbb{R}^{2}}k_{ij}({\bf x-\bf y},t)d{\bf x}\right)f({\bf y})d{\bf y}&&\nonumber\\
=
\left(\int_{\mathbb{R}^{2}}k_{ij}({\bf x},t)d{\bf x}\right)\left(\int_{\mathbb{R}^{2}}f({\bf x})d{\bf x}\right).&&\nonumber\\
&&\label{cdl9}
 \end{eqnarray}
Now, if ${\bf u}_{0}\in L^{1}\times L^{1}$ then the solution of (\ref{cdl3}) satisfies ${\bf u}(\cdot,t)\in L^{1}\times L^{1}, t\geq 0$; the application of (\ref{cdl9}) to (\ref{cdl4}) and hypothesis (H3) imply (\ref{cdl8}).
$\Box$
\begin{remark}
Property (\ref{cdl8}) may be considered as the cross-diffusion version of the average grey-level invariance, when the image is represented by only one real-valued function ( Weickert \cite{Weickert2}). The definition, from ${\bf M}$, of a scalar to play a similar role of average grey-level in this case may however depend on the initial distribution {${\bf u}_{0}=(u_{0},v_{0})^{T}$} of an original real image $f$ in two components. For example, in the case of linear complex diffusion, {${\bf u}_{0}=(f,0)^{T}$} is typically taken and due to the properties of the fundamental solution {(\ref{cdl5})} (Gilboa et al. \cite{GilboaZS2001,GilboaSZ2002,GilboaSZ2004}),
\begin{eqnarray*}
\int_{\mathbb{R}^{2}}h_{R}({\bf x},t)d{\bf x}=1,\quad
\int_{\mathbb{R}^{2}}h_{I}({\bf x},t)d{\bf x}=0,
\end{eqnarray*}
for $t\geq 0$, we have that $I=I_{R}+iI_{I}$ of (\ref{cdl5}) satisfies
\begin{eqnarray*}
m(I_{R}(\cdot,t))=m(f),\quad m(I_{I}(\cdot,t))=0.
\end{eqnarray*}
Then the role of average grey-level might be played by the integral of the real part of the image, that is the first component in the corresponding formulation as a cross-diffusion system. Some other choices of the initial distribution may however {motivate the change of the definition of average grey-level, e.~g.  $m(u)+m(v)$.}
\end{remark}
\begin{remark}
(Flat kernels). A second consequence of mass preserving hypothesis (H3) is that
\begin{eqnarray*}
\lim_{t\rightarrow\infty}k_{ij}(\cdot,t)=0,\quad  i,j=1,2,
\end{eqnarray*}
which means that the kernels are flat as $t\rightarrow \infty$, see Weickert et al. \cite{WeickertII,WeickertII2}.
\end{remark}
\subsubsection{Rotational invariance}
The invariance of the image by rotation is obtained in a similar way to that of the scalar case, see Pauwels et al. \cite{PauwelsVFN1995}.
\begin{lemma}
\label{lem2}
Assume that: 
\begin{itemize}
\item[(H4)]  {For any $t>0$, $ i,j=1,2$, there exists $\kappa_{ij}(\cdot,t)\in L^{1}$ such that $k_{ij}({\bf x},t)=\kappa_{ij}(|{\bf x}|,t)$, ${\bf x} \in \mathbb{R}^2$}
\end{itemize}
{Let $T_{\theta}:\mathbb{R}^{2}\rightarrow\mathbb{R}^{2}$ be a rotation matrix of angle $\theta\in \mathbb{R}$ and for ${\bf u}\in X\times X$ let us define $\mathcal{T}_{\theta}:X\times X\rightarrow X\times X$ as}
\begin{eqnarray*}
(\mathcal{T}_{\theta}{\bf u})({\bf x})={\bf u}(T_{\theta}{\bf x}),\quad {\bf x}\in\mathbb{R}^{2}.
\end{eqnarray*}
Then for any ${\bf u}_{0}\in X\times X$ and $t>0$,
\begin{eqnarray}
\mathcal{K}_{t}(\mathcal{T}_{\theta}{\bf u}_{0})=\mathcal{T}_{\theta}(\mathcal{K}_{t}{\bf u}_{0}).\label{cdl10}
\end{eqnarray}
\end{lemma}

{\em Proof}. The same arguments as those of the scalar case are applied here, since (H4) implies $k_{ij}(T_{\theta}{\bf x},t)=k_{ij}({\bf x},t), t\geq 0, \theta\in \mathbb{R}$ and this leads to (\ref{cdl10}) $\Box$.
\begin{remark}
As in the scalar case, (H4) also implies that for $i,j=1,2, t>0$, {there exists $\widetilde{\kappa}_{ij}=\widetilde{\kappa}_{ij}(\cdot,t)\in L^{1}$ such that
$
\widehat{k}_{ij}({\bf \xi})=\widetilde{\kappa}_{ij}(|{\bf \xi}|), \xi\in\mathbb{R}^{2}.
$ Moreover
} 
\begin{eqnarray*}
\widehat{k}_{ij}({\bf \xi},t)
=\int_{0}^{\infty} \kappa_{ij}(\rho,t)J_{0}(\rho|\xi|)d\rho,\quad t>0,
\end{eqnarray*}
where $J_{0}(z)$ is the zeroth order Bessel function, see Pauwels et al. \cite{PauwelsVFN1995}.
\end{remark}
\subsubsection{Recursivity (semigroup property)}
Here we assume that: 
\begin{itemize}
\item[(H5)] The family of operators $\{\mathcal{K}_{t}, t\geq 0\}$ satisfies the semigroup conditions:
\begin{eqnarray*}
&&\lim_{t\rightarrow 0+}\mathcal{K}_{t}{\bf f}={\bf f},\quad {\bf f}\in X\times X,\\ &&\mathcal{K}_{t+s}=\mathcal{K}_{t}\mathcal{K}_{s}, \quad t,s\geq 0.
\end{eqnarray*}
\end{itemize}
Note that the first assumption in (H5) means that for $i,j=1,2$
\begin{eqnarray*}
\lim_{t\rightarrow 0+}k_{ij}(\cdot,t)=\left\{\begin{matrix}0&\,\,i\neq j\\\delta(\cdot)&\,\,i=j,\end{matrix}\right. 
\end{eqnarray*}
where $\delta(\cdot)$ denotes the Dirac's delta distribution. In terms of the convolution matrices, (H5) reads
\begin{eqnarray}
&&K(\cdot,0)=\delta(\cdot)I,\nonumber\\ 
&&(K(\cdot,t+s)\ast{\bf f})({\bf x})=(K(\cdot,t)\ast K(\cdot,s)\ast{\bf f})({\bf x}),\nonumber\\
&& {\bf x}\in\mathbb{R}^{2}, {\bf u}\in X\times X.\label{cdl11}
\end{eqnarray}
Finally, if we define
\begin{eqnarray}
\widehat{K}(\xi,t)=\begin{pmatrix}\widehat{k}_{11}(\xi,t)&\widehat{k}_{12}(\xi,t)\\
\widehat{k}_{21}(\xi,t)&\widehat{k}_{22}(\xi,t)\end{pmatrix},\quad \xi\in\mathbb{R}^{2}, t\geq 0,\label{cdl12}
\end{eqnarray}
then the second condition in  (\ref{cdl11}) can be written as
\begin{eqnarray}
\widehat{K}(\xi,t+s)=\widehat{K}(\xi,t)\widehat{K}(\xi,s), \quad t,s\geq 0.\label{cdl13}
\end{eqnarray}
On the other hand, the formulation of (\ref{cdl3}) as an initial value-problem of  {PDE} for ${\bf u}$ is related to the existence of an infinitesimal generator of the semigroup $\{\mathcal{K}_{t}, t\geq 0\}$, defined as (Pazy \cite{Pazy1983} and Yosida \cite{Yosida1995})
\begin{eqnarray}
D{\bf f}=\lim_{h\rightarrow 0^{+}}\frac{{\mathcal {K}_{h}}{\bf f}-{\bf f}}{h}, \quad {\bf f}\in X\times X,\label{cdl14}
\end{eqnarray}
and with domain $\mathrm{dom}(D)$ consisting of those functions ${\bf f}\in X\times X$ for which the limit in (\ref{cdl14}) exists. The existence of $D$
is guaranteed under the additional hypothesis of continuity:
\begin{itemize}
\item[(H6)] If $||\cdot ||{_{\ast}}$ is the norm defined in (\ref{cdl1c})  then
\begin{eqnarray*}
\lim_{t\rightarrow 0+}  ||K(\cdot,t)-K(\cdot,0)||{_{\ast}}=0.
\end{eqnarray*}
\end{itemize}
\begin{lemma}
\label{lem3}
Under the hypotheses (H1), (H2), (H5) and (H6) {the function} ${\bf u}(t)={\bf u}(\cdot,t)$ given by (\ref{cdl3}) is a weak solution of the initial value-problem 
\begin{eqnarray}
{\bf u}_{t}(t)&=&D {\bf u}(t),\quad t>0,\nonumber\\
{\bf u}(0)&=&{\bf u}_{0}({\bf x}),\label{cdl15}
\end{eqnarray}
where $D$ is the linear operator (\ref{cdl14}).
\end{lemma}

{\em Proof}. Note that for ${\bf f}\in X\times X$
\begin{eqnarray*}
||K(\cdot,t)\ast {\bf f}-K(\cdot,0)\ast {\bf f}||_{X\times X}
&\leq &\\
||K(\cdot,t)-K(\cdot,0)||||{\bf f}||_{X\times X}&&
\end{eqnarray*}
Therefore, by (H5), (H6) we have that $\{\mathcal{K}_{t}, t\geq 0\}$ is a $C_{0}-$semigroup of bounded linear operators on $X$. This implies that ${\bf u}(t)=\mathcal{K}_{t}{\bf u}_{0}$ is the unique weak solution of (\ref{cdl15}). (If ${\bf u}_{0}\in \mathrm{dom}(D)$ then ${\bf u}(t)$ is {a strong solution}, e.~g. Pazy \cite{Pazy1983}.)
$\Box$

It is worth mentioning the Fourier representation of the matrix operator (\ref{cdl14}). If $D=(D_{ij})_{i,j=1,2}$ denotes the entries of $D$, then the system (\ref{cdl15}) can be written in terms of the Fourier symbols $\widehat{D}(\xi):=\{\widehat{D}_{ij}(\xi)\}_{i,j=1,2}$, as the linear evolution problem
\begin{eqnarray}
\frac{d}{dt}\begin{pmatrix}\widehat{u}({\bf \xi},t)\\ \widehat{v}({\bf \xi},t)\end{pmatrix}&=&\begin{pmatrix}\widehat{D}_{11}({\bf \xi})&\widehat{D}_{12}({\bf \xi})\\\widehat{D}_{21}({\bf \xi})&\widehat{D}_{22}({\bf \xi})\end{pmatrix}\begin{pmatrix}\widehat{u}({\bf \xi},t)\\ \widehat{v}({\bf \xi},t)\end{pmatrix},\label{cdl16}\\
&&{\bf \xi}\in\mathbb{R}^{2}, \quad t>0,\nonumber
\end{eqnarray}
with $\widehat{u}({\bf \xi},0)=\widehat{u}_{0}({\bf \xi}), \widehat{v}({\bf \xi},0)=\widehat{v}_{0}({\bf \xi})$. By taking the Fourier transform in (\ref{cdl4}) and comparing with (\ref{cdl16}) we can write (\ref{cdl12}) in the form
\begin{eqnarray}
\widehat{K}({\bf \xi},t)=e^{t\widehat{D}({\bf \xi})},\quad {\bf \xi}\in\mathbb{R}^{2}, \quad t>0.
\label{cdl17}
\end{eqnarray}
\begin{remark}
\label{rem7}
We also note that the hypothesis for rotational invariance (H4) implies that the matrix $\widehat{D}(\xi)$ only depends on $|\xi|$. Furthermore, (H1) implies that $\widehat{k}_{ij}(\cdot,t)\in C_{0}(\mathbb{R}^{2}), i,j=1,2$. This leads to
\begin{eqnarray*}
\lim_{|{\bf \xi}|\rightarrow\infty}\widehat{D}_{ij}({\bf \xi})\rightarrow -\infty, \quad i,j=1,2,
\end{eqnarray*}
which forces $\widehat{D}(\xi)$ to be negative definite for all $\xi\in\mathbb{R}^{2}$ (that is, the real part of its eigenvalues must be negative for all $\xi\in\mathbb{R}^{2}$).
\end{remark}
\subsubsection{Scale invariance}
The last property to formulate is the scale invariance (\'Alvarez et al. \cite{AlvarezGLM1993}). We define
\begin{eqnarray*}
S_{\lambda}{\bf f}({\bf x})={\bf f}(\lambda{\bf x}),\quad {\bf x}\in\mathbb{R}^{2},\quad {\bf f}\in X\times X,
\end{eqnarray*}
and assume that:
\begin{itemize}
\item[(H7)] For any $\lambda\in\mathbb{R}$ and $t>0$ there is $t^{\prime}= \phi(t)$ such that
\begin{eqnarray*}
S_{\lambda}\mathcal{K}_{t^{\prime}}=\mathcal{K}_{t}S_{\lambda}.
\end{eqnarray*}
\end{itemize}
In Lindeberg \& ter Haar Romeny \cite{LindebergrTH1994} it is argued that, in the context of image processing, the relation between $t$ and the scale represented by the standard deviation $\sigma$ in the Gaussian filtering ($t=\sigma^{2}/2$) can be generalized and assumed to exist from the beginning of the process, by establishing the existence of time ($t$) and scale ($\sigma$) parameters and some connection ($\varphi$) between them. Following this argument, we first introduce a scale parameter $\sigma$, related to the semigroup parameter $t$ by a suitable transformation (to be specified later)
\begin{eqnarray}
t=\varphi(\sigma).\label{cdl18}
\end{eqnarray}
The second condition in (H5), in terms of $\sigma$, reads
\begin{eqnarray}
K(\cdot,\sigma_{1})\ast K(\cdot,\sigma_{2})=K(\cdot,\varphi^{-1}(\varphi(\sigma_{1})+\varphi(\sigma_{2}))),\label{cdl19}
\end{eqnarray}
while the first one implies $\varphi(0)=0$. In order to preserve the qualitative requirement (which is one of the bases of the scale-space theory, see Lindeberg \cite{Lindeberg2009}) that increasing values of the scale parameter {should correspond to representations at coarser scales}, we must assume that $\varphi:\mathbb{R}_{+}\rightarrow\mathbb{R}_{+}$ is monotonically increasing (in particular, invertible). (In Pauwels et al. \cite{PauwelsVFN1995} this $\varphi$ can be identified as $\psi^{-1}$ defined there.)
\subsection{Linear filtering characterization}
\label{sec23}
The characterization of those matrix kernels $K(\cdot,t)$, $t\geq 0$, such that the scale-space representation (\ref{cdl3}) satisfies shift invariance, rotational invariance, recursivity and scale invariance is established in this section. Assume that in terms of the scale $\sigma$, (\ref{cdl3}) is written in the matrix form
\begin{eqnarray*}
{\bf F}(\cdot,\sigma)=K(\cdot,\sigma)\ast {\bf f},
\end{eqnarray*}
where ${\bf f}=(f,g)^{T}, {\bf F}=(F,G)^{T}$;
in Fourier representation this is
\begin{eqnarray}
\widehat{{\bf F}}({\bf \xi},\sigma)=\widehat{K}({\bf \xi},\sigma)\widehat{{\bf f}}({\bf \xi}).\label{cdl20}
\end{eqnarray}
\begin{theorem}
\label{the1} If we assume that $\{\mathcal{K}_{t}, t\geq 0\}$, defined  by (\ref{cdl3}), satisfies (H1)-(H7) then {there exist $p>0$ and a positive definite matrix $d$ such that}
the kernels $K(\cdot,t), t>0$, must have the Fourier representation (\ref{cdl17}) with
\begin{eqnarray}
\widehat{D}(\xi)=-|\xi|^{p}d.\label{cdl21}
\end{eqnarray}
\end{theorem}

{\em Proof}. The proof follows an adaptation of the steps given in Lindeberg \& ter Haar Romeny \cite{LindebergrTH1994}, for the scalar case.
\begin{itemize}
\item[(A)] {\em Dimensional analysis.} The assumption about independence of the scale will allow a simplification of (\ref{cdl20}) by using dimensional analysis, see e.~g. Yarin \cite{Yarin2012}. In this case, taking e.~g. the dimensionless variables $\xi_{1}\sigma, \xi_{2}\sigma, \widehat{f}_{1}\widehat{f}_{2}^{-1}$ (where ${\bf \xi}=(\xi_{1},\xi_{2})^{T}\in\mathbb{R}^{2}$) and applying the Pi-Theorem, {there is a matrix $\widetilde{K}:\mathbb{R}^{2}\rightarrow\mathbb{R}^{2}\times\mathbb{R}^{2}$ with $\widetilde{K}({\bf 0})=I$ (in order to have $\widehat{{\bf F}}({\bf \xi},0)=\widehat{{\bf f}}({\bf \xi})$) such that}
the system (\ref{cdl20}) can be written in the form
\begin{eqnarray*}
\widehat{{\bf F}}({\bf \xi},\sigma)=\widetilde{K}({\bf \xi}\sigma)\widehat{{\bf f}}({\bf \xi}).
\end{eqnarray*}
Furthermore, rotational invariance implies that 
$\widetilde{K}({\bf \xi}\sigma)=\widetilde{K}(|{\bf \xi}\sigma|)$ and therefore
\begin{eqnarray}
\widehat{{\bf F}}({\bf \xi},\sigma)=\widetilde{K}(|{\bf \xi}\sigma|)\widehat{{\bf f}}({\bf \xi}),\label{cdl22}
\end{eqnarray}

\item[(B)] {\em Scale invariance}. According to (\ref{cdl19}) and (\ref{cdl22}), the semigroup condition (H5) is of the form
\begin{eqnarray}
&&\widetilde{K}(|{\bf \xi}\sigma_{1}|)\widetilde{K}(|{\bf \xi}\sigma_{2}|)
=\nonumber\\
&&\widetilde{K}(|{\bf \xi}\varphi^{-1}(\varphi(\sigma_{1})+\varphi(\sigma_{2}))|),\label{cdl23}
\end{eqnarray}
for $\sigma_{1},\sigma_{2}\geq 0$. The same arguments as those in Lindeberg \& ter Haar Romeny \cite{LindebergrTH1994} can be used to show that scale invariance implies that $\varphi$ in (\ref{cdl18}) must be of the form
$$\varphi(\sigma)=C\sigma^{p},$$ for some constant $C>0$ (which can be taken as $C=1$) and $p>0$. (In Pauwels et al. \cite{PauwelsVFN1995}, $p$ is identified as $\alpha$.) Hence, if $H(x^{p})\equiv\widetilde{K}(x)$ then (\ref{cdl23}) reads
\begin{eqnarray*}
H(|{\bf \xi}\sigma_{1}|^{p})H(|{\bf \xi}\sigma_{2}|^{p})&=&\widetilde{K}(|{\bf \xi}\sigma_{1}|)\widetilde{K}(|{\bf \xi}\sigma_{2}|)\\
&=&\widetilde{K}(|{\bf \xi}\varphi^{-1}(\varphi(\sigma_{1})+\varphi(\sigma_{2}))|)\\
&=&\widetilde{K}((|{\bf \xi}\sigma_{1}|^{p}+|{\bf \xi}\sigma_{2}|^{p})^{1/p})\\
&=&H(((|{\bf \xi}\sigma_{1}|^{p}+|{\bf \xi}\sigma_{2}|^{p})^{1/p})^{p})\\
&=&H(|{\bf \xi}\sigma_{1}|^{p}+|{\bf \xi}\sigma_{2}|^{p}),
\end{eqnarray*}
which is identified as the functional equation $$\Psi(\alpha_{1})\Psi(\alpha_{1})=\Psi(\alpha_{1}+\alpha_{2})$$ characterizing the matrix exponential function. Therefore $\widetilde{K}$ must be of the form
\begin{eqnarray*}
\widetilde{K}(|{\bf \xi}\sigma|)=H(|{\bf \xi}\sigma|^{p})=e^{|{\bf \xi}\sigma|^{p}A}, \quad p>0,
\end{eqnarray*}
for some $2\times 2$ real matrix $A$. Now the arguments given in Remark \ref{rem7} show that $A$ must be negative definite or, alternatively
\begin{eqnarray}
\widehat{K}({\bf \xi},\sigma)=\widetilde{K}(|{\bf \xi}\sigma|)=e^{-|{\bf \xi}\sigma|^{p}d},\quad \xi\in\mathbb{R}^{2},\label{cdl24}
\end{eqnarray}
{where $d$ is a $2\times 2$ positive definite matrix}. Writing (\ref{cdl24}) in terms of the original scale $t$ leads to the representation (\ref{cdl17}) with $\widehat{D}(\xi)$ given by (\ref{cdl21}). $\Box$
\end{itemize}
\begin{remark}
The form (\ref{cdl21}) corresponds to specific forms of the infinitesimal generator $D$. Note that if ${\bf f}\in X\times X$ then
\begin{eqnarray*}
\widehat{\left(\frac{\mathcal{K}_{h}{\bf f}-{\bf f}}{h}\right)}(\xi)=\left(\frac{e^{-h|{\bf \xi}|^{p}d}-I}{h}\right)\widehat{{\bf f}}({\bf \xi}),
\end{eqnarray*}
and formally
\begin{eqnarray*}
\frac{e^{-h|{\bf \xi}|^{p}d}-I}{h}&=&\sum_{j=1}^{\infty}
\frac{(-1)^{j}h^{j-1}|{\bf \xi}|^{jp}}{j!}d^{j}.
\end{eqnarray*}
Thus
\begin{eqnarray*}
\lim_{h\rightarrow 0^{+}}\frac{e^{-h|{\bf \xi}|^{p}d}-I}{h}=-|{\bf \xi}|^{p}d,
\end{eqnarray*}
and the limit is the Fourier symbol of the operator, \cite{Pazy1983}
\begin{eqnarray}
D{\bf f}=-(-\Delta)^{p/2}d{\bf f},\label{cdl25}
\end{eqnarray}
with $\Delta$ standing for the Laplace operator and where $(-\Delta)^{p/2}$ is multiplying each entry of $d$, cf. Pauwels et al. \cite{PauwelsVFN1995}.

The explicit formula (\ref{cdl25}) can be used to discuss the additional scale-space property of locality (Weickert et al. \cite{WeickertII,WeickertII2}). A semigroup of operators $T_{t}, t\geq 0$, satisfies the locality condition if for all smooth ${\bf f}, {\bf g}$ in its domain and all ${\bf x}\in \mathbb{R}^{2}$
\begin{eqnarray*}
(T_{t}{\bf f}-T_{t}{\bf g})({\bf x})=o(t),\quad t\rightarrow 0^{+},
\end{eqnarray*}
 whenever the derivatives of ${\bf f}$ and ${\bf g}$ of any nonnegative order are identical. Mathematically (\'Alvarez et al. \cite{AlvarezGLM1993}, Pauwels et al. \cite{PauwelsVFN1995}), the locality condition implies that the corresponding  infinitesimal generator is a local differential operator, which in the case of (\ref{cdl25}) means that $p/2$ must be integer, extending the result obtained in Pauwels et al. \cite{PauwelsVFN1995} for the scalar case to the cross-diffusion framework.
\end{remark}
\begin{remark}
Note that when $d$ has complex conjugate eigenvalues,  Lemma \ref{lem_app} in the appendix shows that there is a basis in $\mathbb{R}^{2}$ such that $d$ is similar to a matrix of the form
\begin{eqnarray}
\begin{pmatrix}\nu&-\mu\\\mu&\nu\end{pmatrix},\label{cdl28}
\end{eqnarray}
with $\nu>0, \mu\neq 0$ and where the eigenvalues of $d$ are $\nu\pm i\mu$. Therefore, linear complex diffusion  {corresponds} to the case of (\ref{cdl25}) with $p=2$ and a matrix $d$ of the form (\ref{cdl26}). The complex diffusion coefficient  {is} given by $c=\nu+i\nu$ or $c=\nu-i\mu$. Formula (\ref{cdl25}) shows that this linear complex diffusion can be generalized by using other values of $p$.
\end{remark}
\begin{remark}
The nature of the semigroup $\{\mathcal{K}_{t}, t\geq 0\}$ can be analyzed from the spectrum and regularity of the infinitesimal generator (\ref{cdl25}) (Pazy \cite{Pazy1983}). In this sense, the following result holds.
\begin{theorem}
\label{the2}
Let $k\geq 0$. {The operator $D$ in (\ref{cdl25}) with domain $\mathrm{dom}(D)=H^{k+1}\times H^{k+1}$ is the infinitesimal generator of a $C_{0}$-semigroup $\{\mathcal{K}_{t}, t\geq 0\}$ and there exists $M>0$ such that the induced norm satisfies
$||\mathcal{K}_{t}||\leq M$}. Furthermore, if $d$ is of one of the three reduced forms 
\begin{eqnarray}
d_{1}&=&\begin{pmatrix}\lambda_{+}&0\\0&\lambda_{-}\end{pmatrix},\quad \lambda_{+}\geq \lambda_{-}>0,\nonumber\\
d_{2}&=&\begin{pmatrix}\alpha&\beta\\0&\alpha\end{pmatrix},\quad \alpha\geq \beta>0,\nonumber\\
d_{3}&=&\begin{pmatrix}\nu&-\nu\\\mu&\nu\end{pmatrix},\quad \nu>0, \mu\neq 0,
\end{eqnarray}
 then $M\leq 1$.
\end{theorem}

{\em Proof}. We first prove that for each reduced form $d_{j}, j=1,2,3$ in (\ref{cdl28}), $\{\mathcal{K}_{t}, t\geq 0\}$ is a $C_{0}$-semigroup of contractions. Consider the eigenvalue problem for $D$:
\begin{eqnarray}
(\lambda I-D){\bf u}={\bf f},\label{cdl29}
\end{eqnarray}
where $I$ is the $2\times 2$ identity matrix, ${\bf u}=(u,v)^{T}, {\bf f}=(f,g)^{T}$. Assume that in (\ref{cdl21}) $d=d_{1}$. In terms of the Fourier transform, (\ref{cdl29}) reads
\begin{eqnarray*}
(\lambda+|\xi|^{p}\lambda_{+})\widehat{u}({\bf \xi})&=&\widehat{f}({\bf \xi}),\nonumber\\
(\lambda+|\xi|^{p}\lambda_{-})\widehat{v}({\bf \xi})&=&\widehat{g}({\bf \xi}),\quad \xi\in \mathbb{R}^{2}.
\end{eqnarray*}
Then, since $\lambda_{+}>\lambda_{-}>0$, for any $\lambda>0$ we have
\begin{eqnarray}
\frac{1}{|\lambda+|\xi|^{p}\lambda_{\pm}|}\leq \frac{1}{\lambda},\quad \xi\in \mathbb{R}^{2}.\label{cdl30}
\end{eqnarray}
When $d=d_{2}$, the Fourier system associated to (\ref{cdl29}) is now of the form
\begin{eqnarray*}
(\lambda+|\xi|^{p}\alpha)\widehat{u}({\bf \xi})+|\xi|^{p}\beta\widehat{v}({\bf \xi})&=&\widehat{f}({\bf \xi}),\nonumber\\
(\lambda+|\xi|^{p}\alpha)\widehat{v}({\bf \xi})&=&\widehat{g}({\bf \xi}),\quad \xi\in \mathbb{R}^{2}.
\end{eqnarray*}
Now, since $0<\beta<\alpha$ then for any $\lambda>0$,
\begin{eqnarray}
\frac{|\xi|^{p}\beta}{|\lambda+|\xi|^{p}\alpha|}\leq 1,\quad 
\frac{1}{|\lambda+|\xi|^{p}\alpha|}\leq \frac{1}{\lambda},\quad \xi\in \mathbb{R}^{2}.\label{cdl31}
\end{eqnarray}
Finally, when $d=d_{3}$, the Fourier representation of (\ref{cdl29}) has the form
\begin{eqnarray}
(\lambda+\nu|{\bf \xi}|^{p})\widehat{u}({\bf \xi})-\mu|{\bf \xi}|^{p}\widehat{v}({\bf \xi})&=&\widehat{f}({\bf \xi}),\nonumber\\
\mu|{\bf \xi}|^{p}\widehat{u}({\bf \xi})+(\lambda+\nu|{\bf \xi}|^{p})\widehat{v}({\bf \xi})&=&\widehat{g}({\bf \xi}).\label{cdl32}
\end{eqnarray}
Inverting (\ref{cdl32}) leads to
\begin{eqnarray*}
\widehat{u}({\bf \xi})&=&\frac{(\lambda+\nu|{\bf \xi}|^{p})}{m({\bf \xi})}\widehat{f}({\bf \xi})+\frac{\mu|{\bf \xi}|^{p}}{m({\bf \xi})}\widehat{g}({\bf \xi}),\\
\widehat{v}({\bf \xi})&=&-\frac{\mu|{\bf \xi}|^{p}}{m({\bf \xi})}\widehat{f}({\bf \xi})+\frac{(\lambda+\nu|{\bf \xi}|^{p})}{m({\bf \xi})}\widehat{g}({\bf \xi}),\\
m({\bf \xi})&=&(\lambda+\nu|{\bf \xi}|^{p})^{2}+(\mu|{\bf \xi}|^{p})^{2}.
\end{eqnarray*}
Note now that since $\nu>0$, for any $\lambda>0$ we have
\begin{eqnarray}
\frac{(\lambda+\nu|{\bf \xi}|^{p})}{m({\bf \xi})}\leq \frac{1}{\lambda+\nu|{\bf \xi}|^{p}}\leq \frac{1}{\lambda},\label{cdl33}
\end{eqnarray} and also, since
\begin{eqnarray*}
|\lambda+\nu|{\bf \xi}|^{p}||\mu|{\bf \xi}|^{p}|\leq \frac{m({\bf \xi})}{2},
\end{eqnarray*} 
then
\begin{eqnarray}
\frac{|\mu|{\bf \xi}|^{p}|}{m({\bf \xi})}\leq \frac{1}{2(\lambda+\nu|{\bf \xi}|^{p})}\leq \frac{1}{2\lambda}.\label{cdl34}
\end{eqnarray}
Finally,  the application of Hille-Yosida theorem to each case, using the corresponding estimates (\ref{cdl30}), (\ref{cdl31}) and (\ref{cdl32})-(\ref{cdl34}), proves the second part of the theorem. For the general case, we note that, using the eigenvalues of $d$, there is a nonsingular matrix $P$ such that $\Lambda=PdP^{-1}$ is of one of the three forms in (\ref{cdl28}) (see Lemma \ref{lem_app}). Then the theorem follows by using $M=||P||||P^{-1}||$.
$\Box$
\end{remark}
\subsection{Generalized small theta approximation}
\label{sec24}
One of the arguments to consider complex diffusion as an alternative for image processing is the so-called small theta approximation (Gilboa et al. \cite{GilboaSZ2004}). This means that for small values of the imaginary part of the complex diffusion coefficient, the corresponding imaginary part of the solution of the evolutionary diffusion problem behaves, in the limit, as a scaled smoothed Gaussian derivative of the initial signal. This idea can also be discussed in the context of cross-diffusion systems (\ref{cdl3}), where $D$ is the infinitesimal generator (\ref{cdl18}), that is
\begin{eqnarray}
{\bf u}({\bf x},t)=e^{tD}{\bf u}_{0}({\bf x})=e^{-t(-\Delta)^{p/2}d}{\bf u}_{0}({\bf x}).\label{cdl35}
\end{eqnarray}
By using the notation introduced in Lemma \ref{lem_app}, we have the following result.
\begin{theorem}
\label{the3}
Define the operator $A=-(-\Delta)^{p/2}, p>0$ and let $d$ be a positive definite matrix with eigenvalues and parameters given by (\ref{cdl27}). Assume that $d$ satisfies one of the cases (i), (iii) or (iv) in Lemma \ref{lem_app}. Let $f\in X$ be a real-valued function. If ${\bf u}({\bf x},t)=(u({\bf x},t),v({\bf x},t))^{T}$ then (\ref{cdl35}) satisfies:
\begin{itemize}
\item[(C1)] If $|s|\rightarrow 0$ and ${\bf u}_{0}=(f,0)^{T}$ then
\begin{eqnarray*}
&&\lim_{d_{21}\rightarrow 0}{u({\bf x},t)}=e^{\frac{q}{2} tA}f({\bf x}),\\ 
&&\lim_{d_{21}\rightarrow 0}\frac{v({\bf x},t)}{d_{21}}=tA\left(e^{\frac{q}{2} tA}f({\bf x})\right).
\end{eqnarray*}
\item[(C2)] If $|s|\rightarrow 0$ and ${\bf u}_{0}=(0,f)^{T}$ then
\begin{eqnarray*}
&& \lim_{d_{12}\rightarrow 0}\frac{u({\bf x},t)}{d_{12}}=tA\left(e^{\frac{q}{2} tA}f({\bf x})\right),\\
&&\lim_{d_{12}\rightarrow 0}{v({\bf x},t)}=e^{\frac{q}{2} tA}f({\bf x}),
\end{eqnarray*}
\end{itemize}
where $q=\mathrm{tr}(d)$ is the trace of $d$.
\end{theorem}

{\em Proof}. We can write (\ref{cdl35}) in the form
\begin{eqnarray}
{\bf u}({\bf x},t)=Pe^{tA\Lambda}P^{-1}{\bf u}_{0}({\bf x}),\label{cdl36}
\end{eqnarray}
where $P$ and $\Lambda$ are the corresponding matrices specified in Lemma \ref{lem_app} in each case. Specifically and after some tedious but straightforward calculations, we have:
\begin{itemize}
\item In the case (i)
\begin{eqnarray}
u({\bf x},t)&=&e^{\frac{q}{2} tA}\left(\cosh(t\mu A)-\frac{r}{2\mu}\sinh(t\mu A)\right)u_{0}({\bf x})\nonumber\\
&&-
e^{\frac{q}{2} tA}\left(\frac{d_{12}}{\mu}\sin(t\mu A)\right)v_{0}({\bf x}),\nonumber\\
v({\bf x},t)&=&e^{\frac{q}{2} tA}\left(\cosh(t\mu A)+\frac{r}{2\mu}\sinh(t\mu A)\right)v_{0}({\bf x})\nonumber\\
&&+
e^{\frac{q}{2} tA}\left(\frac{d_{21}}{\mu}\sinh(t\mu A)\right)u_{0}({\bf x}),\label{cdl37}
\end{eqnarray}
\item In the case (iii)
\begin{eqnarray}
u({\bf x},t)&=&e^{\frac{q}{2} tA}\left((1-\frac{rt}{2}A)u_{0}({\bf x})+d_{12}tAv_{0}({\bf x})\right),\nonumber\\
v({\bf x},t)&=&e^{\frac{q}{2} tA}\left((1+\frac{rt}{2}A)v_{0}({\bf x})+d_{21}tAu_{0}({\bf x})\right),\nonumber\\\label{cdl38}
\end{eqnarray}
\item In the case (iv)
\begin{eqnarray}
u({\bf x},t)&=&e^{\frac{q}{2} tA}\left(\cos(t\mu A)-\frac{r}{2\mu}\sin(t\mu A)\right)u_{0}({\bf x})\nonumber\\
&&-
e^{\frac{q}{2} tA}\left(\frac{d_{12}}{\mu}\sin(t\mu A)\right)v_{0}({\bf x}),\nonumber\\
v({\bf x},t)&=&e^{\frac{q}{2} tA}\left(\cos(t\mu A)+\frac{r}{2\mu}\sin(t\mu A)\right)v_{0}({\bf x})\nonumber\\
&&+
e^{\frac{q}{2} tA}\left(\frac{d_{21}}{\mu}\sin(t\mu A)\right)u_{0}({\bf x}) ,\label{cdl39}
\end{eqnarray}
\end{itemize}
where ${\bf u}_{0}({\bf x})=(u_{0}({\bf x}),v_{0}({\bf x}))^{T}$ and the cosine, sine, hyperbolic cosine and hyperbolic sine of the operator $A$ are defined in the standard way from the exponential, {see e.~g. Yosida \cite{Yosida1995}}. By using the approximations as $z\rightarrow 0$,
\begin{eqnarray*}
\cos(z)\approx 1,\quad \sin(z)\approx z,\quad \cosh(z)\approx 1,\quad \sinh(z)\approx z,
\end{eqnarray*}
and the corresponding limits in (\ref{cdl37})-(\ref{cdl39}) then (C1) and (C2) hold.
$\Box$

Theorem \ref{the3} can be considered as a generalization of the small theta approximation property of linear complex diffusion. Under the conditions specified in the theorem, one of the components behaves as the operator $A$ applied to a smoothed version of the original image $f$, scaled by $t$ and with the smoothing effect determined by $A$ and the trace of $d$. Actually, formulas (\ref{cdl37})-(\ref{cdl39}) can be used to extend the result to other initial distributions ${\bf u}_{0}$. Finally, note that if $d$ is similar to a matrix of the case (ii) in Lemma \ref{lem_app} then
\begin{eqnarray*}
u({\bf x},t)=e^{\frac{q}{2} tA}u_{0}({\bf x}),\quad
v({\bf x},t)=e^{\frac{q}{2} tA}v_{0}({\bf x}),
\end{eqnarray*}
and this property does not hold. This case could be considered as a Gaussian smoothing for both components.
\section{Numerical experiments}
\label{sec3}
This section is devoted to illustrate numerically the behaviour of linear filters of cross-diffusion. More specifically, the numerical experiments presented here, involving one- and two-dimensional signals, will concern the influence, in a filtering problem with processes (\ref{cdl3}) and kernels satisfying (\ref{cdl21}), of the following elements:
\begin{itemize}
\item The matrix $d$: According to the analysis of Section \ref{sec2}, the choice of $d$ plays a role in the generalized small theta approximation and the experiments will also attempt to discuss if the influence is somehow extended to the quality of {filtering}. Three types of matrices will be taken in the experiments, covering the different form of the eigenvalues (see Lemma \ref{lem_app}).
\item The initial distribution ${\bf u}_{0}$: Besides the generalized small theta approximation, the choice of ${\bf u}_{0}$ is also relevant to define the average grey value.
\item The parameter $p$: Here the purpose is to explore numerically if locality affects the {filtering} process, either in its quality or computational performance.
\end{itemize}
As mentioned in the Introduction, these experiments are a first approach to the numerical study of the behaviour of the filters. All the computations are made with Fourier techniques (Canuto et al. \cite{CHQZ}). More presicsely, for the case of one-dimensional signals, an interval $(-L,L)$ with large $L$ is defined and discretized by Fourier collocation points $x_{j}=-L+jh, j=0,\ldots,N,$ with stepsize $h>0$ and the signal is represented by the corresponding trigonometric interpolating polynomial with the coefficients computed by the discrete Fourier Transform (DFT) of the signal values at the collocation points. For experiments with images, the implementation follows the corresponding Fourier techniques in two dimensions, with discretized intervals $(-L_{x},L_{x})$ $\times$ $(-L_{y},L_{y})$, being $L_{x}, L_{y}$ large, by Fourier collocation points $(x_{j},y_{k})$, with  $x_{j}=-L_{x}+jh_{x}, j=0,\ldots,N_{x},$ $y_{k}=-L_{y}+kh_{y}, k=0,\ldots,N_{y}$, and the image is represented by the trigonometric interpolating polynomial at the collocation points, computed with the corresponding two-dimensional version of the DFT. In both cases, from the Fourier representation, the convolution (\ref{cdl3}) is implemented in the Fourier space by using (\ref{cdl15}).

The experiments can be divided in two groups. The first one concerns the evolution of a clean signal. It illustrates properties like the generalized small theta approximation and the effect of locality. The second group deals with image filtering problems and the experiments are performed with the following strategy: from an original image $s$ we add some noise to generate an initial noisy signal. (We have made experiments with several types of noise, but only those with Gaussian type will be shown here.) From this noisy signal an initial distribution is defined and then the restored image, given by (\ref{cdl3}), is monitored at several times $t_{n}, n=0,1,\ldots$, in order to estimate the quality of restoration. This quality has been measured by using different standard metrics, namely:
\begin{itemize}
\item Signal-to-Noise-Ratio (SNR):
\begin{eqnarray}
SNR(s,u^{n})=10\log_{10}\left(\frac{{\rm var}(u^{n})}{{\rm var}(s-u^{n})}\right).\label{snr}
\end{eqnarray}
\item Peak Signal-to-Noise-Ratio (PSNR):
 \begin{eqnarray}
PSNR(s,u^{n})=10\log_{10}\left(\frac{l^{2}}{||s-u^{n}||^{2}}\right).\label{psnr}
\end{eqnarray}
\end{itemize}
In all the cases, $u^{n}$ stands for the first component of the restored image at time $t_{n}$, $l$ is the length of the vectors in one-dimensional signals and $l=255$ for two-dimensional signals, $||\cdot ||$ stands for the Euclidean norm in 1D and the Frobenius norm in 2D and ${\rm var}(x)$ is the variance of the vector $x$ (or the matrix disposed as a vector) 
According to the formulas, the larger the parameter values the better the filtering is. Other metrics, like the Root-Mean-Square-Error or the correlation coefficient, have been used in the experiments, although only the results corresponding to (\ref{snr}), (\ref{psnr}) will be shown here.
\subsection{Choice of matrix $d$}
\label{sec31}
\subsubsection{Experiments in 1D}
\label{sec311}
\begin{figure}[htbp]
\centering
\subfigure[]
{\includegraphics[width=0.5\textwidth]{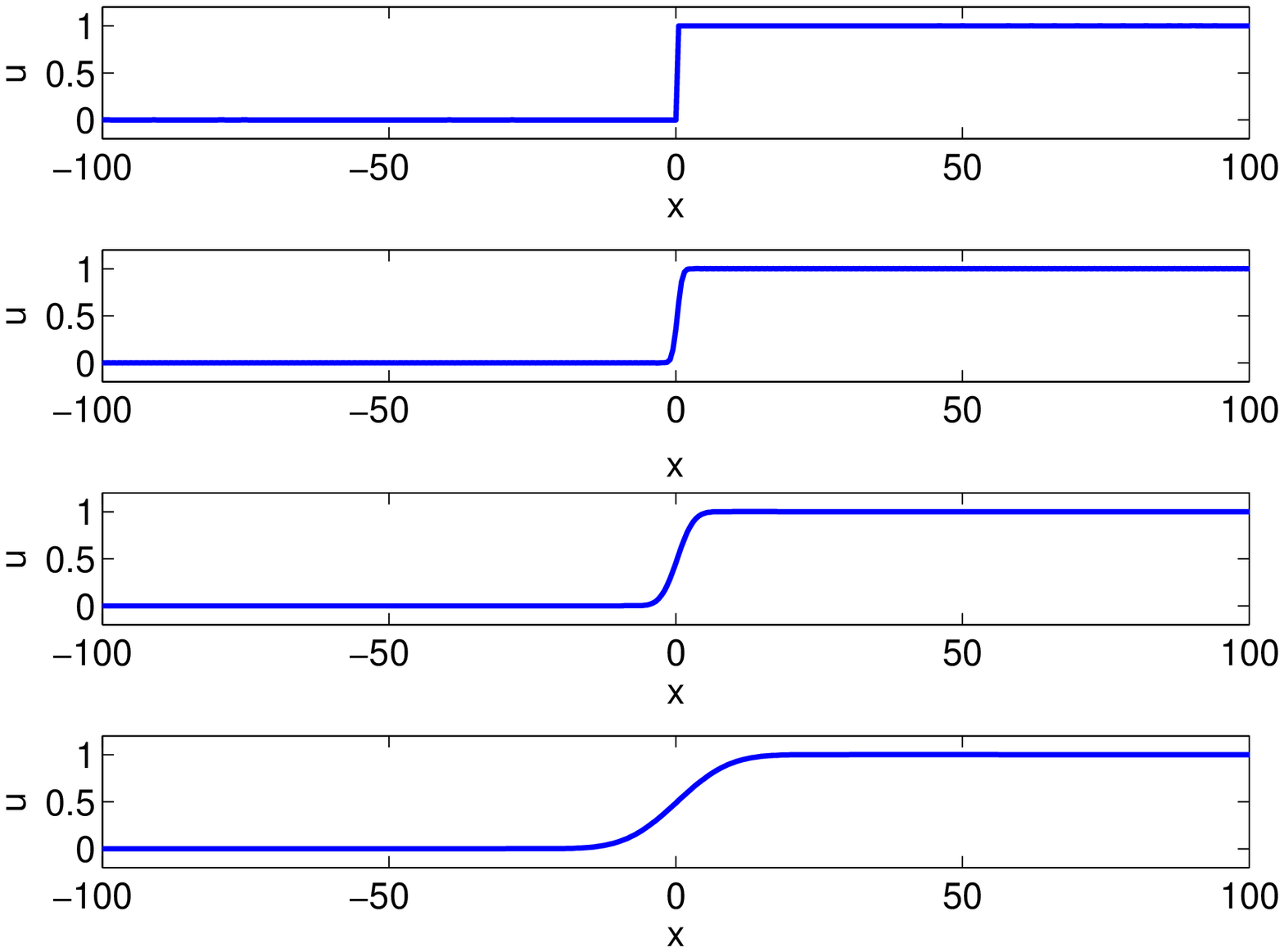}}
\subfigure[]
{\includegraphics[width=0.5\textwidth]{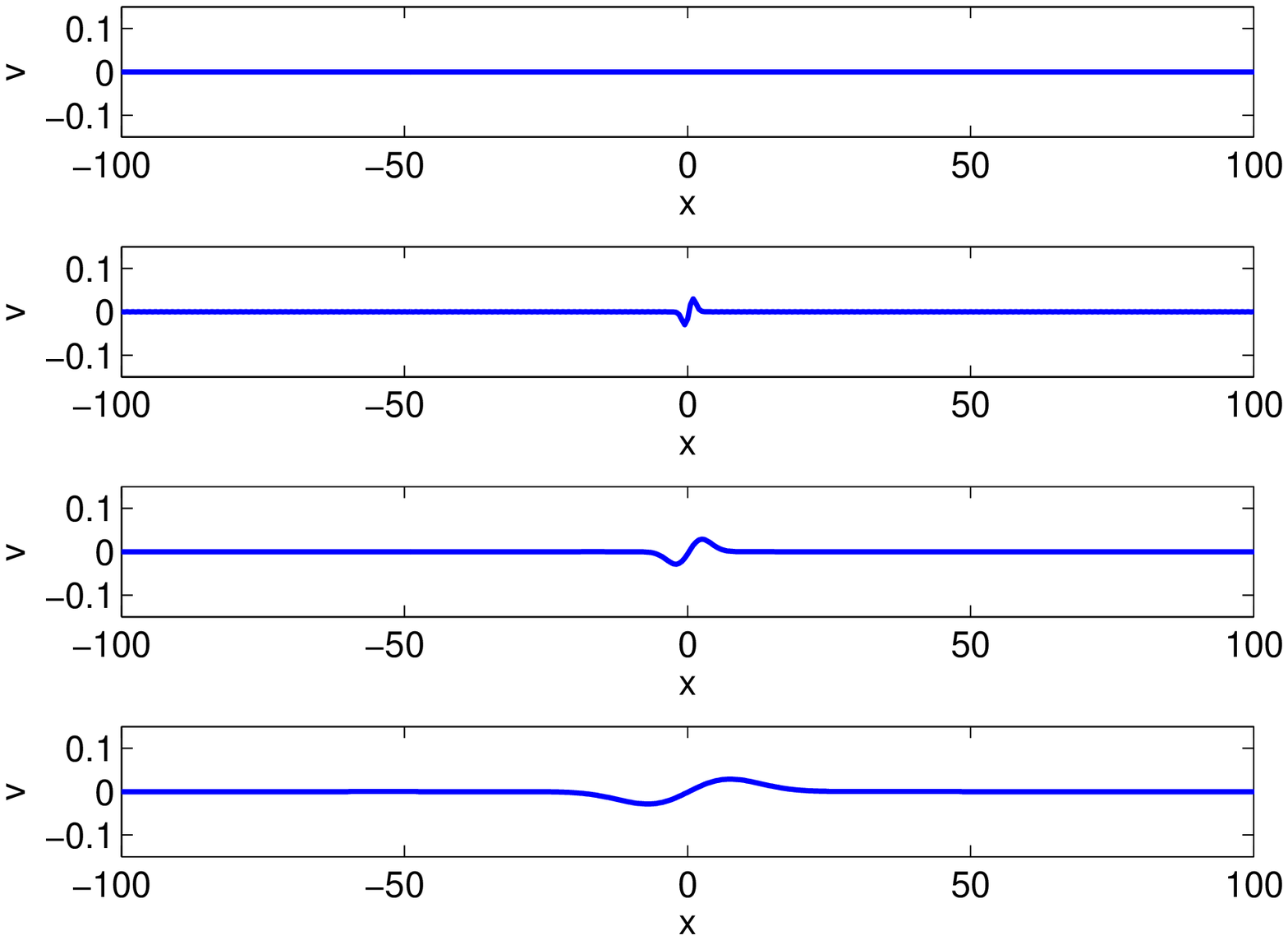}}
\caption{Cross-diffusion with $p=2$ and for $d_{11}=1,d_{22}=1.1, d_{12}=0.1,d_{21}=-0.25$. Profiles of  (a) $u$ and (b) $v$ at times $t=0,0.25,2.5,25$.}
\label{fig_I_2}
\end{figure}

\begin{figure}[htbp]
\centering
\subfigure[]
{\includegraphics[width=0.5\textwidth]{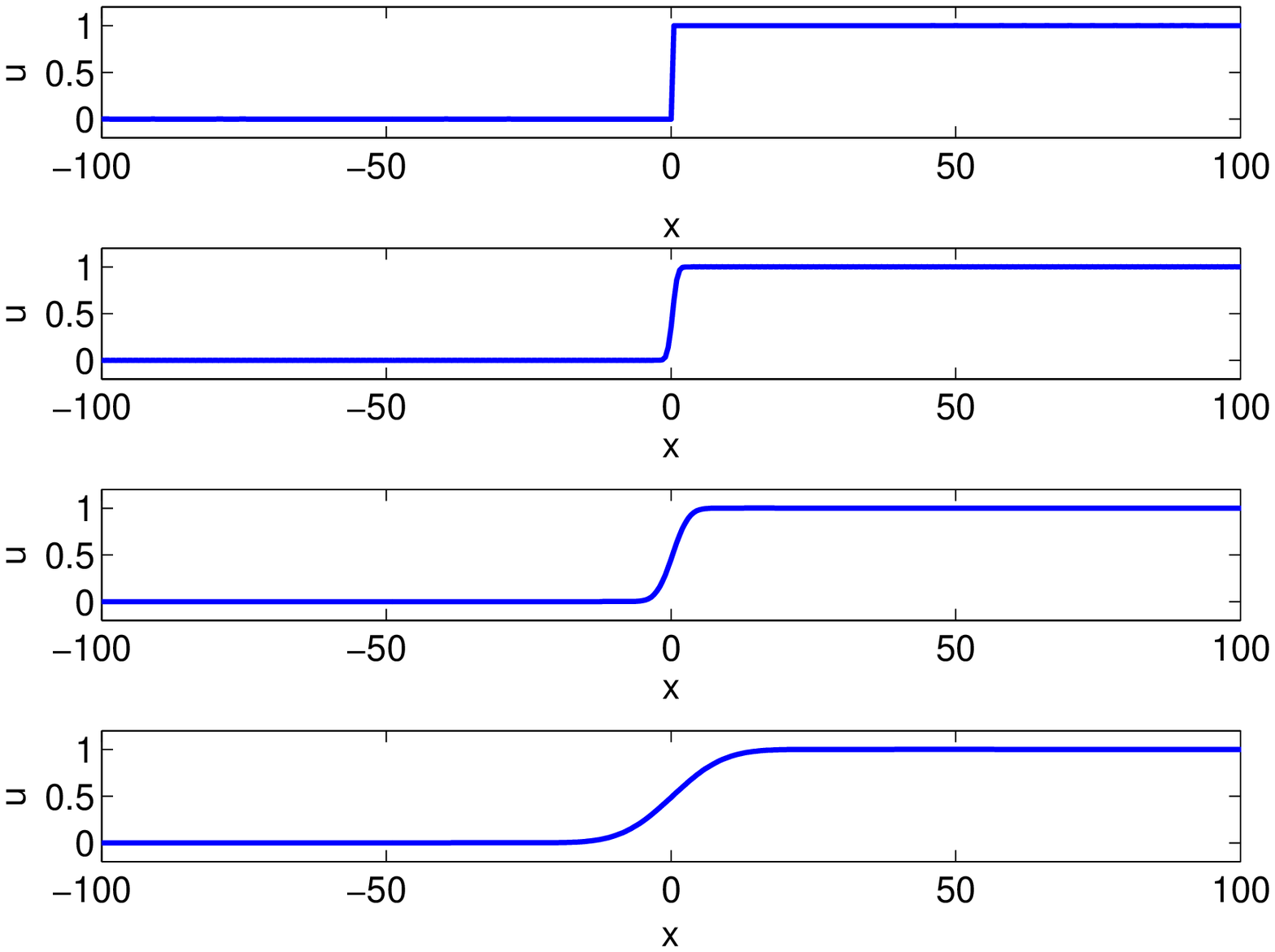}}
\subfigure[]
{\includegraphics[width=0.5\textwidth]{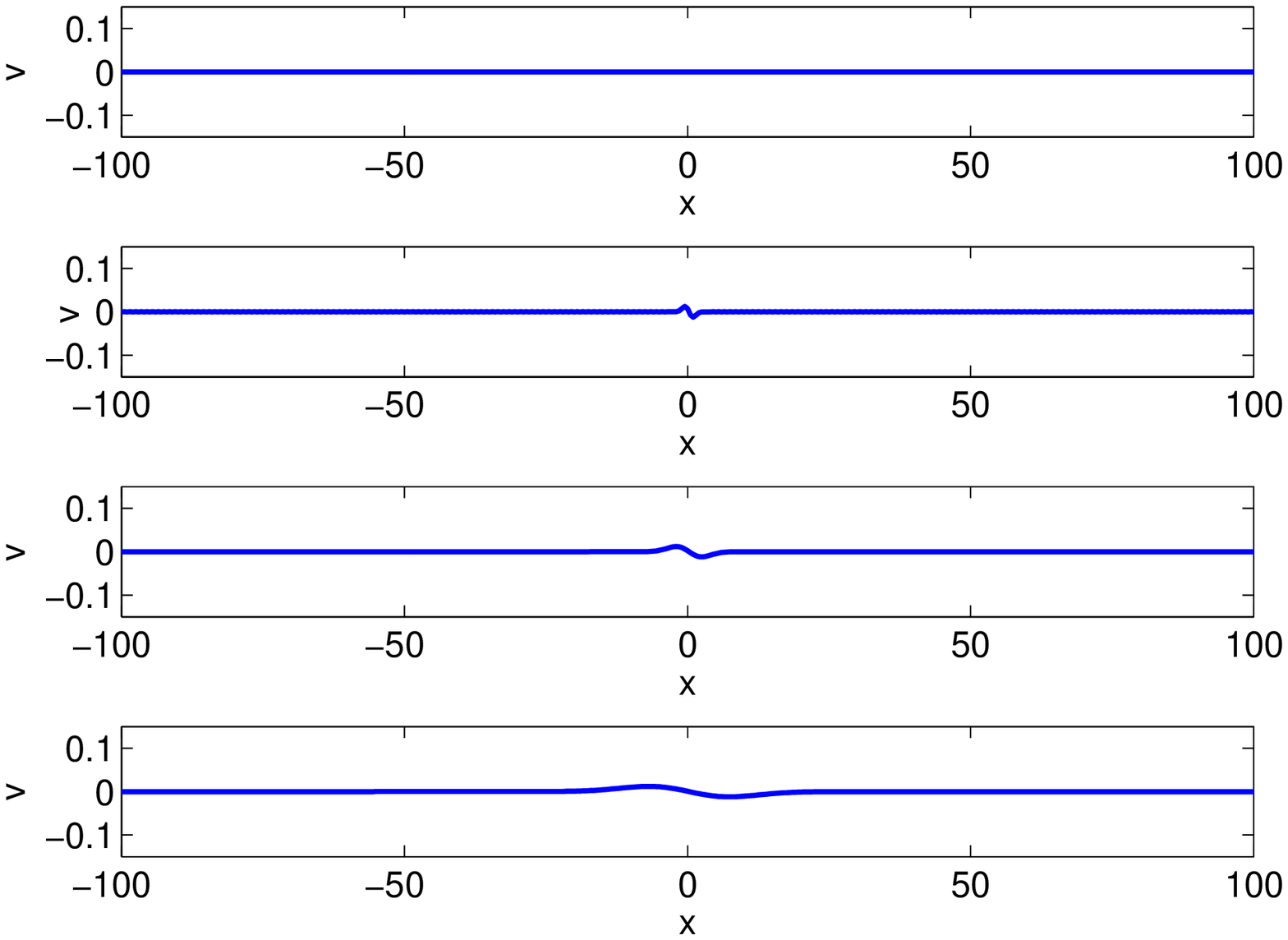}}
\caption{Cross-diffusion with $p=2$ and for $d_{11}=d_{22}=1, d_{12}=-0.1,d_{21}=0.1$. Profiles of  (a) $u$ and (b) $v$ at times $t=0,0.25,2.5,25$.}
\label{fig_I_3}
\end{figure}

{A unit function $f$ is first considered}. Taking ${\bf u}_{0}=(f,0)^{T}$, the evolution (\ref{cdl3}) with kernels satisfying Theorem \ref{the1} and $p=2$ is monitored at several times. Experiments with three types of matrices $d$ (covering the cases $s>0, s=0, s<0$, where $s$ is given by (\ref{cdl27})) and for each type, with different values (according to the size of $|s|$) were performed. They are represented in Figures  \ref{fig_I_1}-\ref{fig_I_3} and suggest two first conclusions:

\begin{figure}[htbp]
\centering
\subfigure[]
{\includegraphics[width=0.5\textwidth]{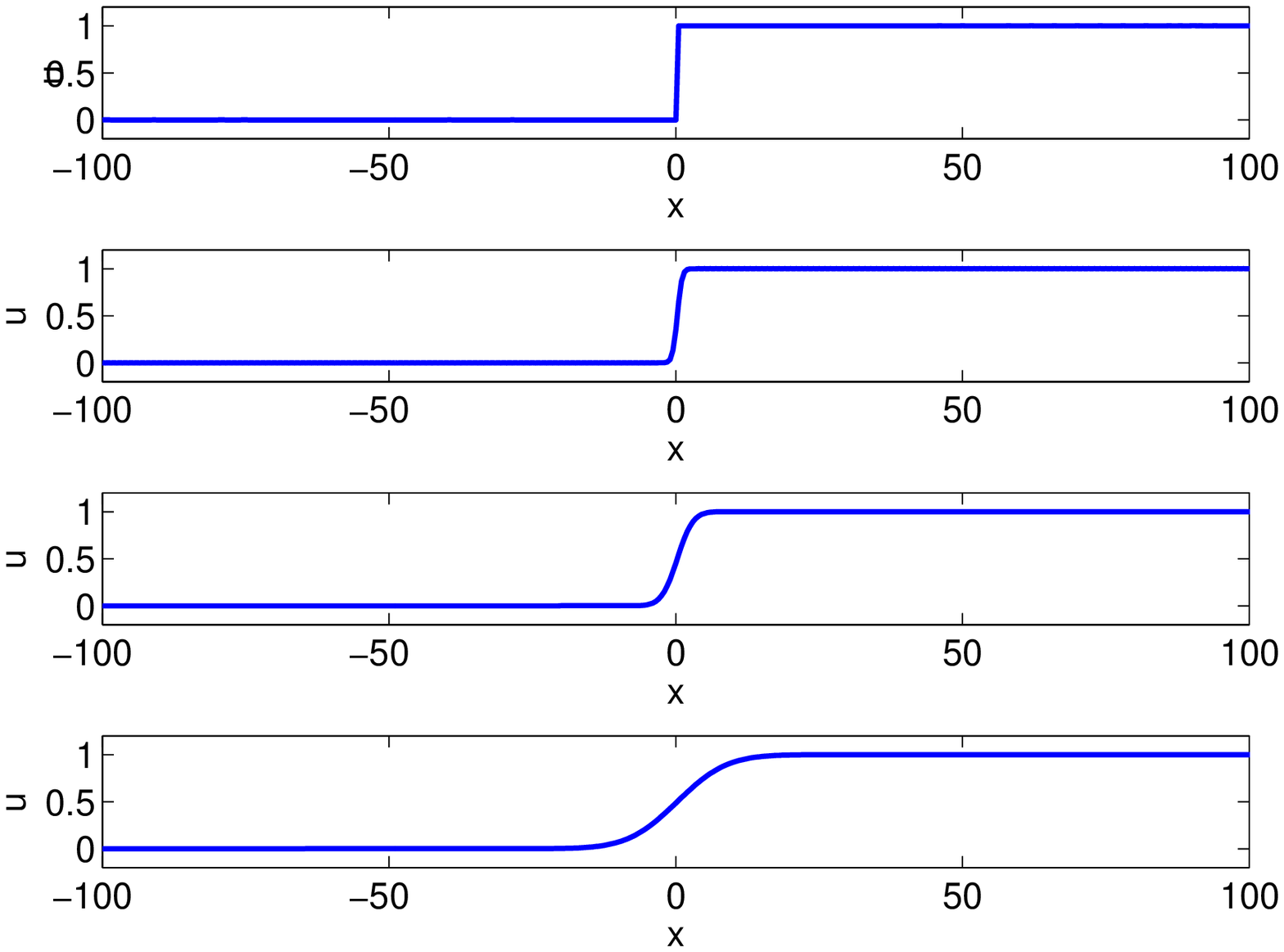}}
\subfigure[]
{\includegraphics[width=0.5\textwidth]{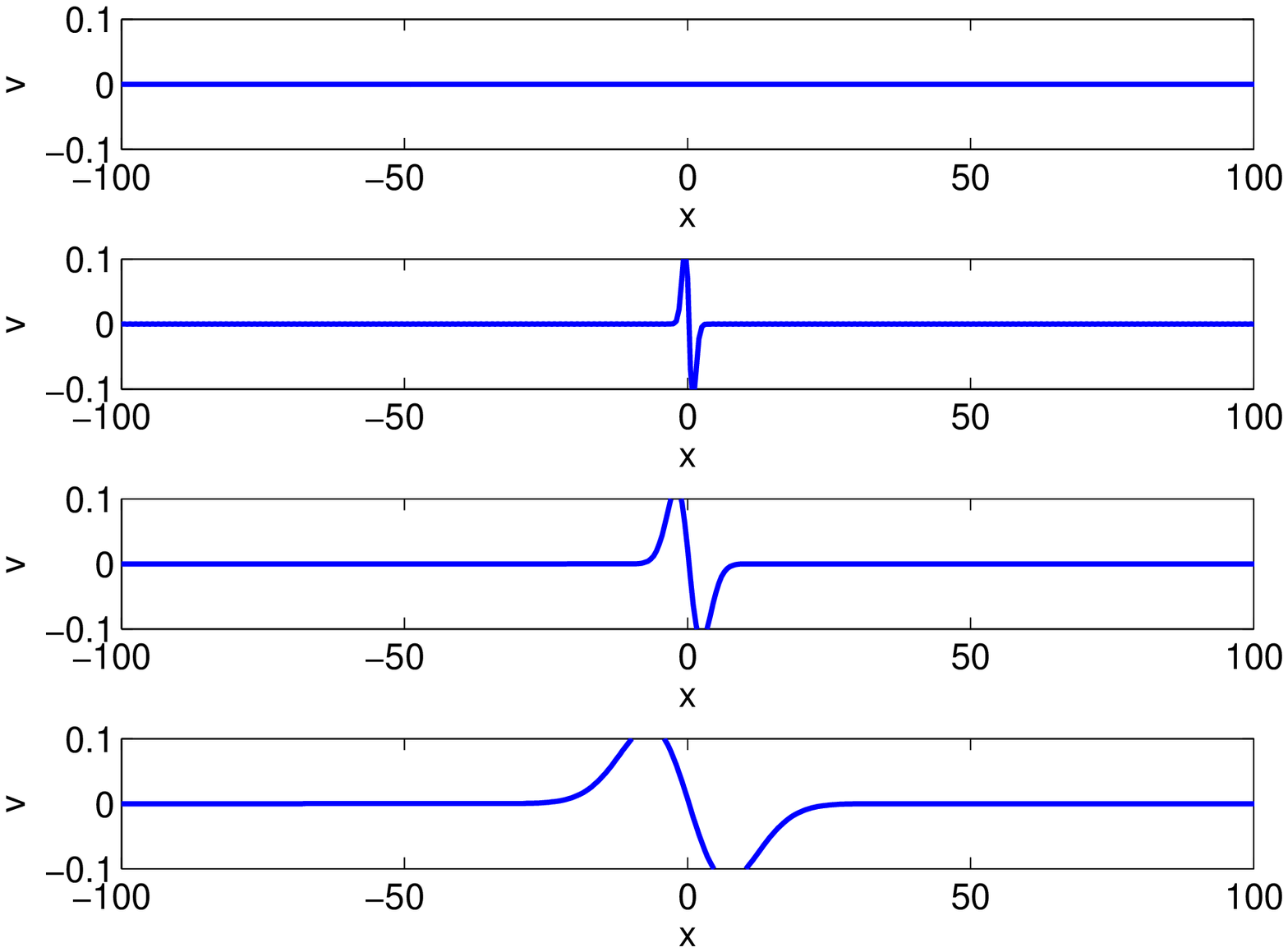}}
\caption{Cross-diffusion with $p=2$ and for $d_{11}=1,d_{22}=1.1, d_{12}=0.1,d_{21}=1$. Profiles of  (a) $u$ and (b) $v$ at times $t=0,0.25,2.5,25$.}
\label{fig_I_1}
\end{figure}

\begin{enumerate}
\item The first component is affected by a smoothing effect. {We observed that the regularization} is stronger as $d_{11}, d_{22}$ (which, by positive definiteness of $d$, must be positive) and $|d_{12}|, |d_{21}|$ grow.
\item Except in the case $d_{12}=d_{21}=0$ (which may be associated to real Gaussian smoothing) the second component develops a sort of small-amplitude Gaussian derivative-type monopulse. Again, the height of the amplitude depends on how large (in absolute value) the elements of $d$ are, with the larger the parameters the taller the wave is. In particular, this property illustrates the effect of the small theta approximation in complex diffusion (Figure \ref{fig_I_3}) and in more general cross-diffusion models (Figures \ref{fig_I_1} and \ref{fig_I_2}). 
\end{enumerate}
\subsubsection{Experiments in 2D}
\label{sec312}
\begin{figure}[htbp]
\centering
\subfigure[]
{\includegraphics[width=0.5\textwidth]{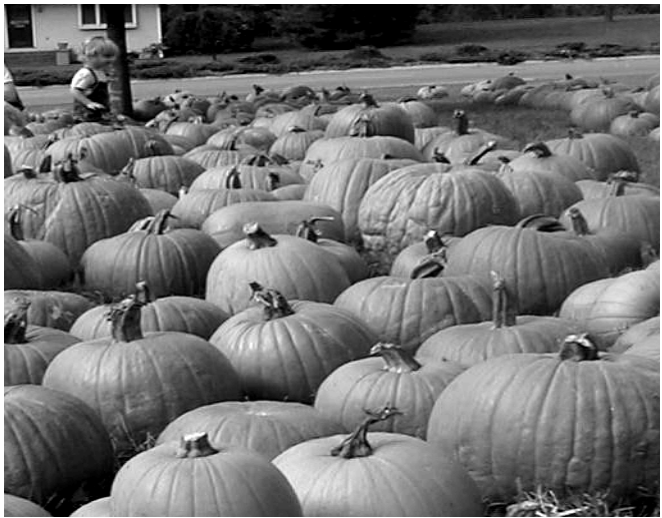}}
\subfigure[]
{\includegraphics[width=0.5\textwidth]{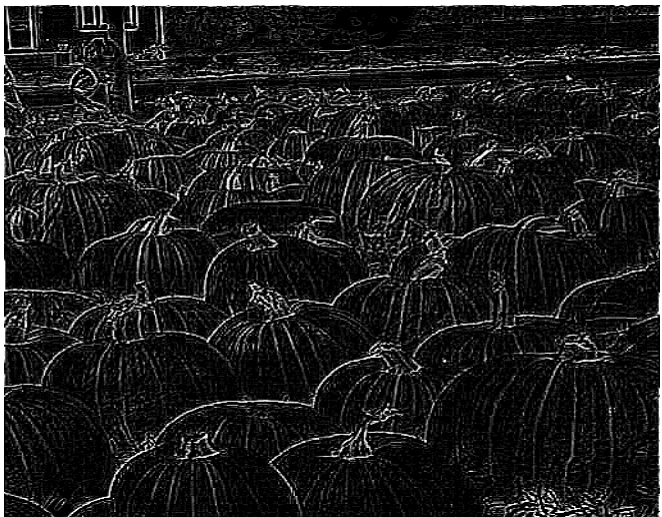}}
\subfigure[]
{\includegraphics[width=0.5\textwidth]{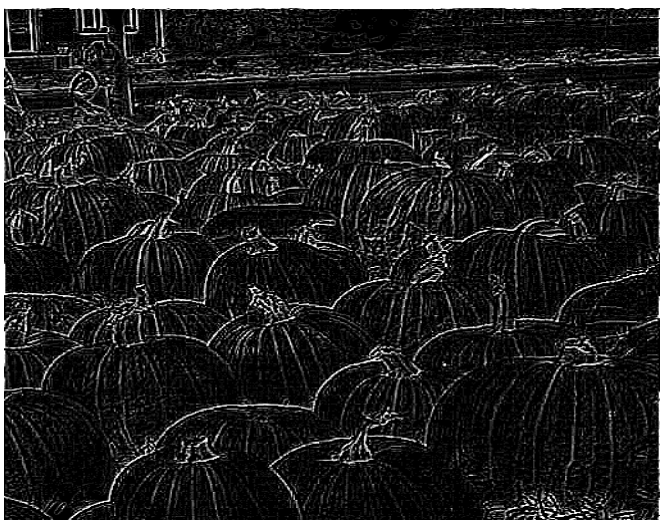}}
\caption{Generalized small theta approximation. (a) Original image $f$. (b) and (c) Second component of the solution of (\ref{cdl3}) at $t=0.1$ with $p=2, {\bf u}_{0}=(f,0)^{T}$ and (b) $d=d_{1}$, (c) $d=d_{2}$, see (\ref{cdl40}).}
\label{figR_A}
\end{figure}

\begin{figure}[htbp]
\centering
\subfigure[]
{\includegraphics[width=0.5\textwidth]{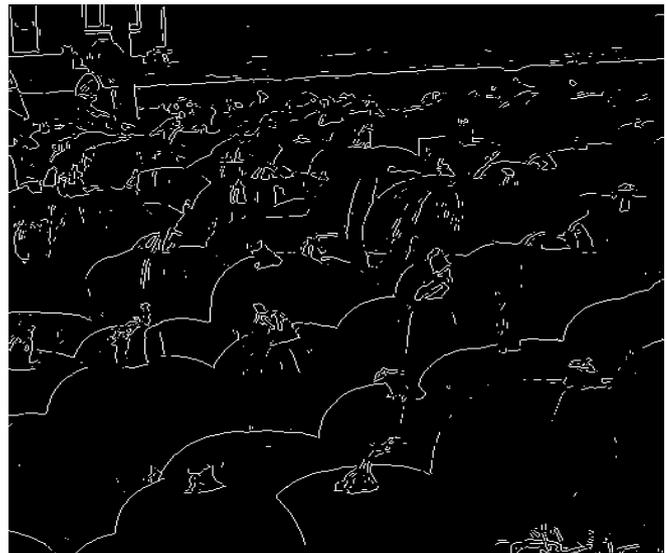}}
\subfigure[]
{\includegraphics[width=0.5\textwidth]{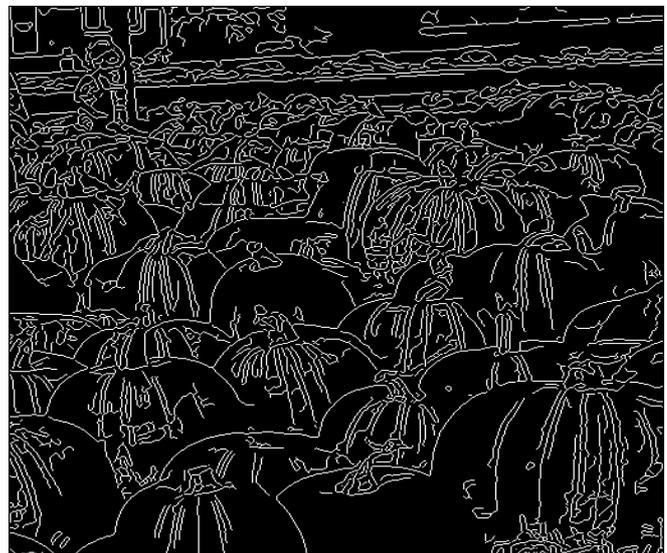}}
\caption{Edge detection from $f$ in Figure \ref{figR_A}(a) provided by: (a) the Prewitt method (Prewitt \cite{Prewitt1970}) and (b) the Canny method (Canny, \cite{Canny1986}).}
\label{figR_B}
\end{figure}

\begin{figure}[htbp]
\centering
\subfigure[]
{\includegraphics[width=0.5\textwidth]{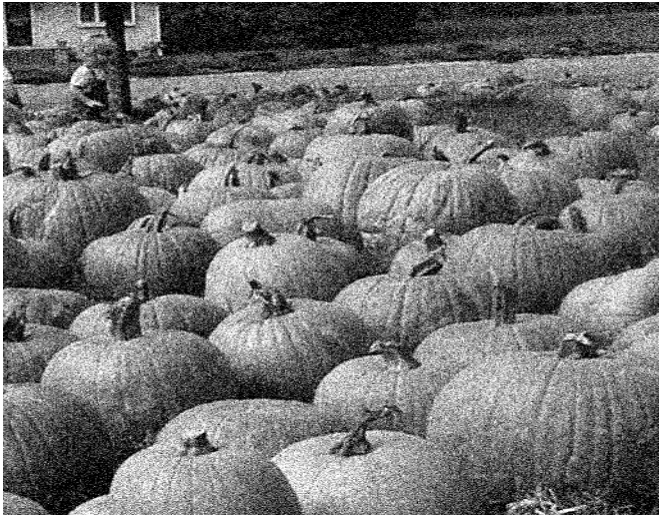}}
\subfigure[]
{\includegraphics[width=0.5\textwidth]{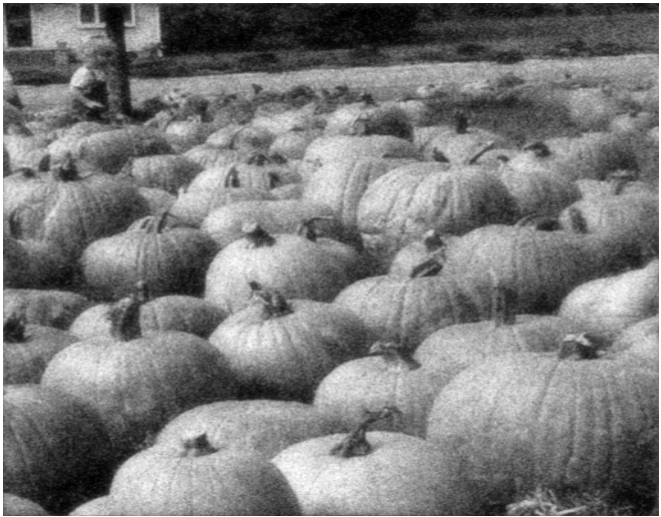}}
\subfigure[]
{\includegraphics[width=0.5\textwidth]{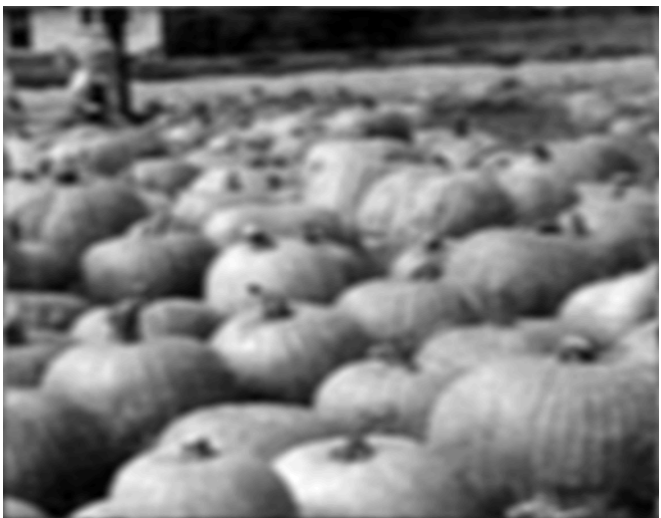}}
\caption{Image filtering problem. (a) Original noisy image $f$. (b) and (c) First component of the solution of (\ref{cdl3}) at $t=15$ with $p=2, {\bf u}_{0}=(f,0)^{T}$ and (b) $d=d_{1}$, (c) $d=d_{2}$, see (\ref{cdl41}).}
\label{figR_C}
\end{figure}

\begin{figure}[htbp]
\centering
\subfigure[]
{\includegraphics[width=0.5\textwidth]{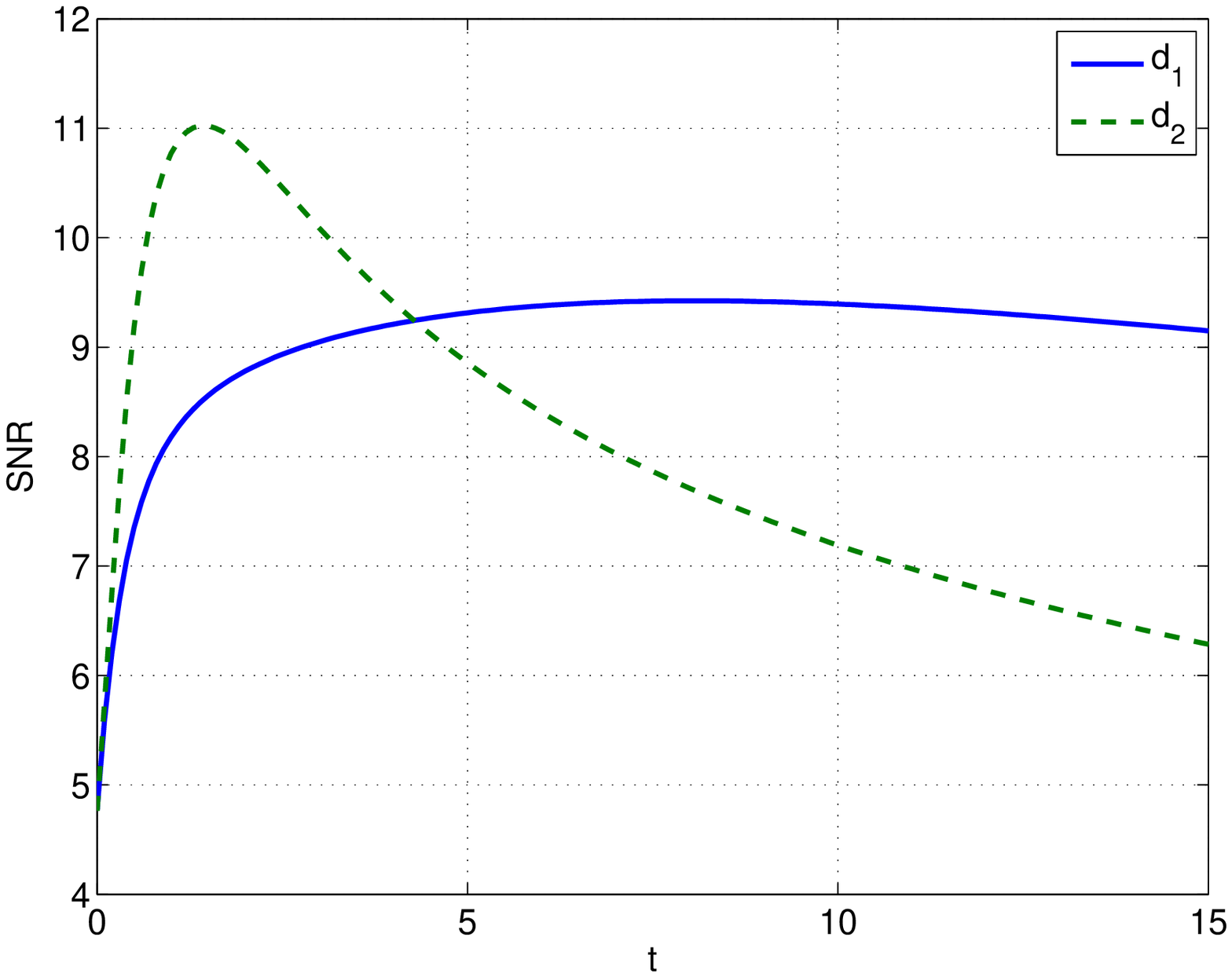}}
\subfigure[]
{\includegraphics[width=0.5\textwidth]{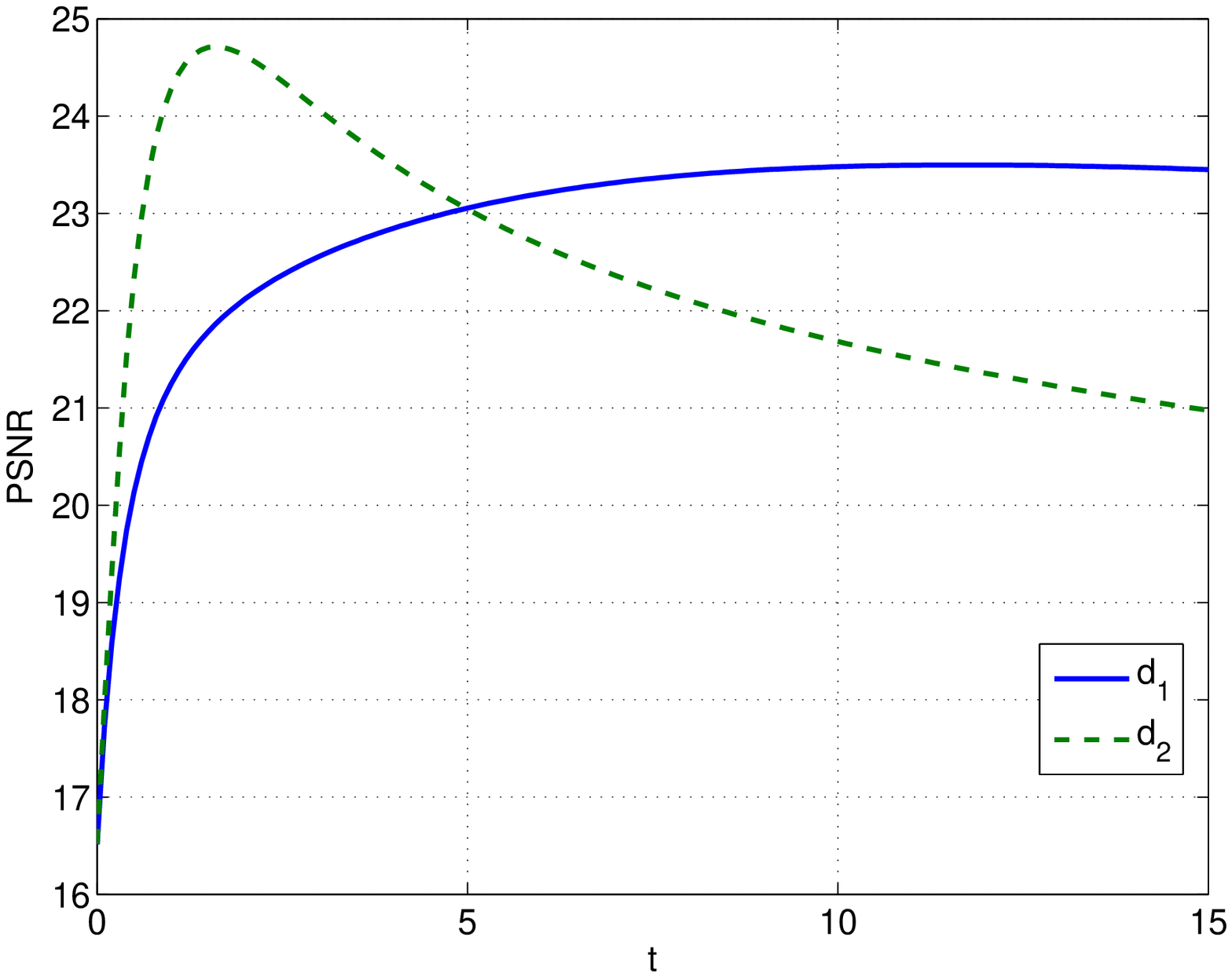}}
\caption{(a) $SNR$ vs. $t$ and (b) $PSNR$ vs. $t$ for a filter (\ref{cdl3}) with $p=2, {\bf u}_{0}=(f,0)^{T}$, $f$ in Figure \ref{figR_C}(a) and {for both  $d=d_{1}$ and  $d=d_{2}$}, see (\ref{cdl41}).}
\label{figR_D}
\end{figure}

The illustration of the influence of the matrix $d$ in 2D signals is focused on mainly two points: the small theta approximation and the behaviour of the filtering evolution of (\ref{cdl3}) with respect to the blurring effect, the detection of the edges and the treatment of textures of the image. From the experiments performed, we observe the following:
\begin{itemize}
\item The property of generalized small theta approximation is illustrated in Figure \ref{figR_A}, which corresponds to apply (\ref{cdl3}) at $t=0.1$ with $p=2$ and matrices
\begin{eqnarray}
d_{1}=\begin{pmatrix}1&10^{-5}\\1.99&1\end{pmatrix},\quad
d_{2}=\begin{pmatrix}1&-10^{-5}\\1.99&1\end{pmatrix},\label{cdl40}
\end{eqnarray}
corresponding to the cases $s>0$ and $s<0$, respectively, with $|s|$ small. The initial condition is ${\bf u}_{0}=(f,0)^{T}$ where $f$ is the original image displayed in Figure \ref{figR_A}(a). According to Theorem \ref{the3}, the entry $d_{21}$ of $d$ may be used as a natural scale for the second component of the solution of (\ref{cdl3}), which is shown in Figures \ref{figR_A}(b) and (c) for $d=d_{1}$ and $d=d_{2}$ respectively. In this small theta approximation property the processes displayed (and those that we performed with matrices for which $s=0$, not shown here) do not have relevant differences: the second component is affected by a slight blurring. A comparison with some standard methods for edge detection (Canny, \cite{Canny1986}, Prewitt \cite{Prewitt1970}), is given in Figure \ref{figR_B}.
\end{itemize}
\begin{itemize}
\item
When the image is affected by some noise, the behaviour of the filtering process with (\ref{cdl3}) can show some differences depending on the choice of the matrix $d$. The first one concerns the blurring effect. Our numerical experiments suggest a better evolution of the filtering for matrices in the case (i) of Lemma \ref{lem_app} than that of matrices in cases (ii)-(iv). (The last one includes the linear complex diffusion.) This is illustrated by Figure \ref{figR_C}, that corresponds to the evolution of (\ref{cdl3}) from a noisy image $f$ (Figure \ref{figR_C}(a)) with $p=2, {\bf u}_{0}=(f,0)^{T}$ and matrices $d=d_{1}$ (Figure \ref{figR_C}(b)) and $d=d_{2}$ (Figure \ref{figR_C}(c)), now of the form
\begin{eqnarray}
d_{1}=\begin{pmatrix}1&0.9\\1&1\end{pmatrix},\quad
d_{2}=\begin{pmatrix}1&-0.9\\1&1\end{pmatrix},\label{cdl41}
\end{eqnarray}
at time $t=15$. The experiment shows that filtering with $d_{1}$ delays the blurring effect with respect to the behaviour observed in the case of $d_{2}$. Similar experiments suggest that using matrices $d$ for which the parameter $s=(d_{22}-d_{11})^{2}+4d_{12}d_{21}$ is positive and moderately large (always conditioned to the satisfaction of positive definite character) improves the filtering in this way. This is also confirmed by Figure \ref{figR_D}, where the evolution of the corresponding SNR and PSNR indexes (\ref{snr}) and (\ref{psnr}) for (\ref{cdl3}) with $d_{1}$ and $d_{2}$ are shown, respectively.
\end{itemize}

\begin{figure}[htbp]
\centering
\subfigure[]
{\includegraphics[width=0.5\textwidth]{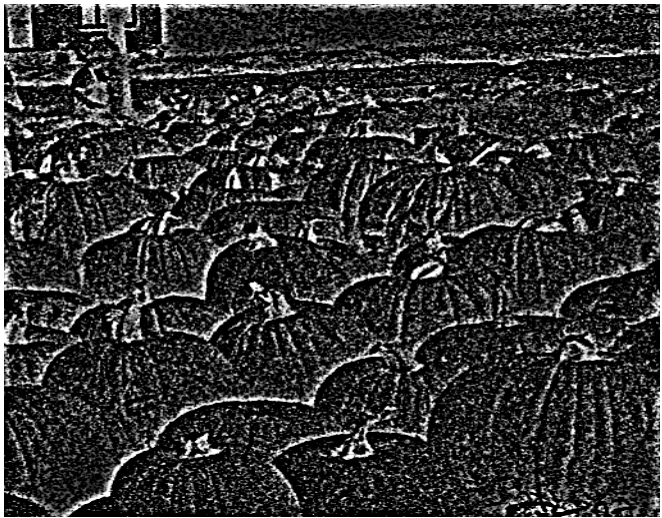}}
\subfigure[]
{\includegraphics[width=0.5\textwidth]{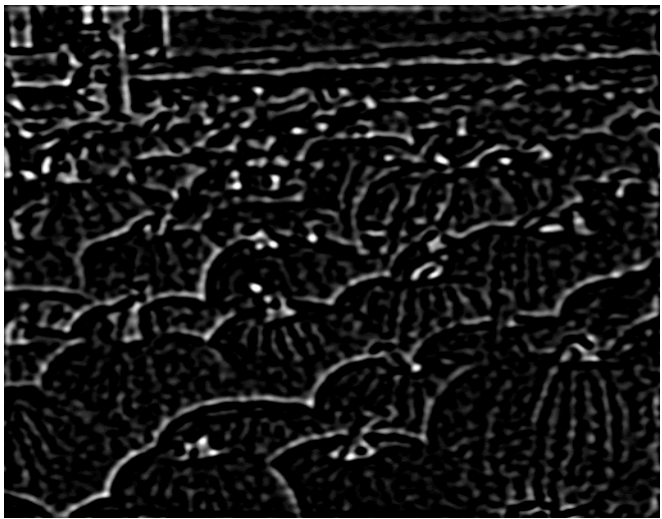}}
\caption{Image filtering problem. Second component of the solution of (\ref{cdl3}) at $t=15$ with $p=2, {\bf u}_{0}=(f,0)^{T}$, $f$ in Figure \ref{figR_C}(a) and: (a) $d=d_{1}$, (b) $d=d_{2}$, see (\ref{cdl41}).}
\label{figR_E}
\end{figure}

\begin{itemize}
\item
The delay of the blurring effect has also influence in other features of the image. The first one is the edge detection by using the second component of (\ref{cdl3}), as observed in Figure \ref{figR_E}. The better performance of the process with $d_{1}$ gives a less blurred detection of the edges than the one given by $d_{2}$.
On the other hand, the delay of the blurring effect may improve the identification of the textures of the image. Using nonlinear models, Lorenz et al. \cite{LorenzBZ} suggested a good behaviour of cross-diffusion with respect to the textures. Numerical experiments in this sense were performed here and they are illustrated in Table \ref{tav_A}. This shows the entropy as a measure of texture. The entropy was computed as
\begin{eqnarray}
En=\sum_{i,j}c_{ij}\log_{2}{c_{ij}},\label{cdl42}
 \end{eqnarray} 
 where $c_{ij}$ stands for the entries of the corresponding grey-level occurrence matrix (  Gonzalez et al. \cite{GonzalezWE2003}).
\end{itemize}
 \begin{table}
\caption{Entropy (\ref{cdl42}) for the image in Figure \ref{figR_A}(a) and the corresponding images obtained with (\ref{cdl3}) at $t=10$ with $p=2, {\bf u}_{0}=(f,0)^{T}$, {being $f$ the corresponding original image, for both  $d=d_{1}$ and  $d=d_{2}$}, see (\ref{cdl41}).}\label{tav_A}
\begin{center}
\begin{tabular}{c|c|c|c}
\hline
&Original&$d_{1}$&$d_{2}$\\
\hline
Figure \ref{figR_A}(a) &$7.5477$&$7.5750$&$7.5200$\\
\hline
\end{tabular}
\end{center}
\end{table}


\subsection{Choice of the parameter $p$}
\label{sec32}
\subsubsection{Experiments in 1D}
\label{sec321}
\begin{figure}[htbp]
\centering
\subfigure[]
{\includegraphics[width=0.5\textwidth]{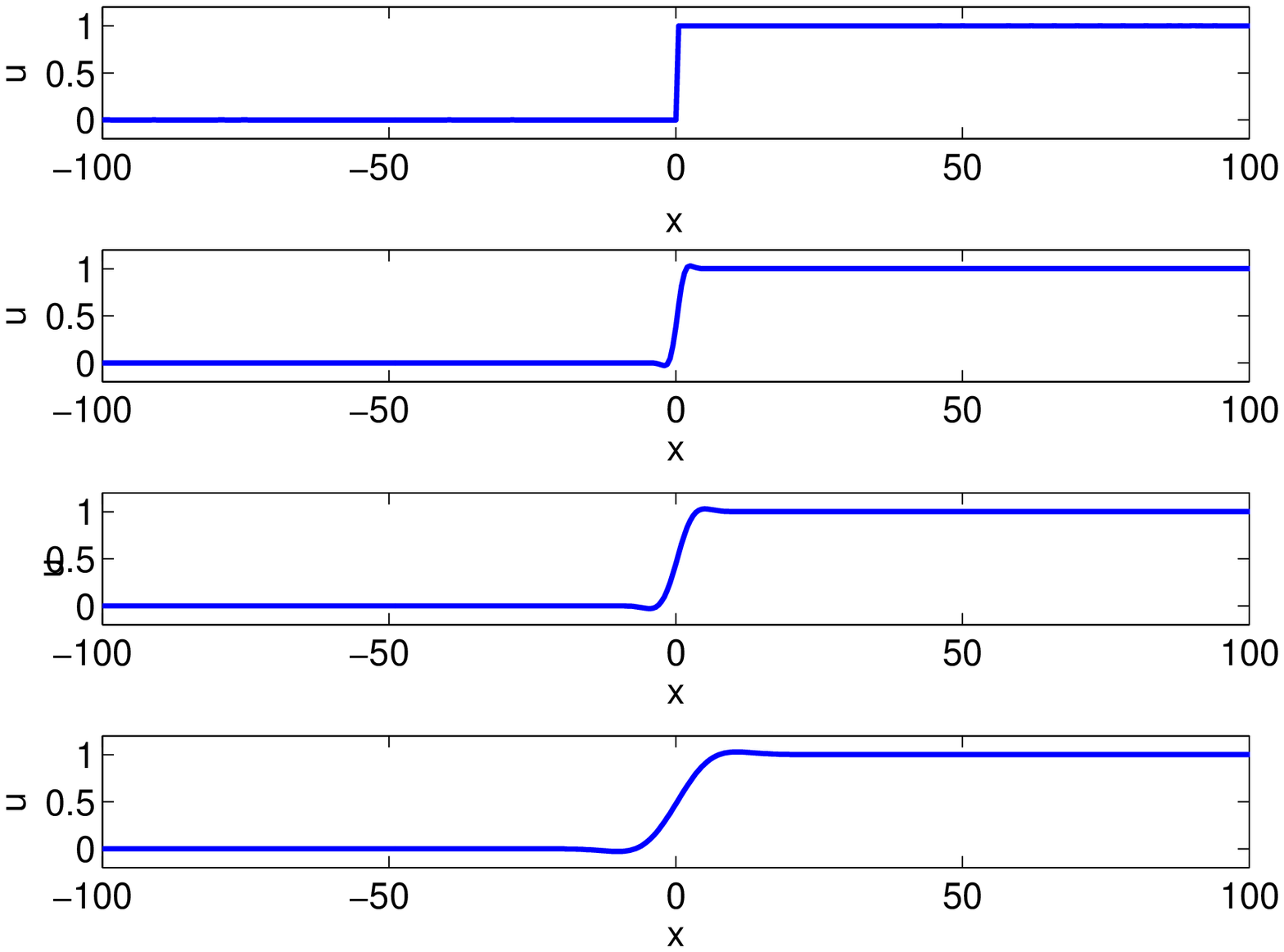}}
\subfigure[]
{\includegraphics[width=0.5\textwidth]{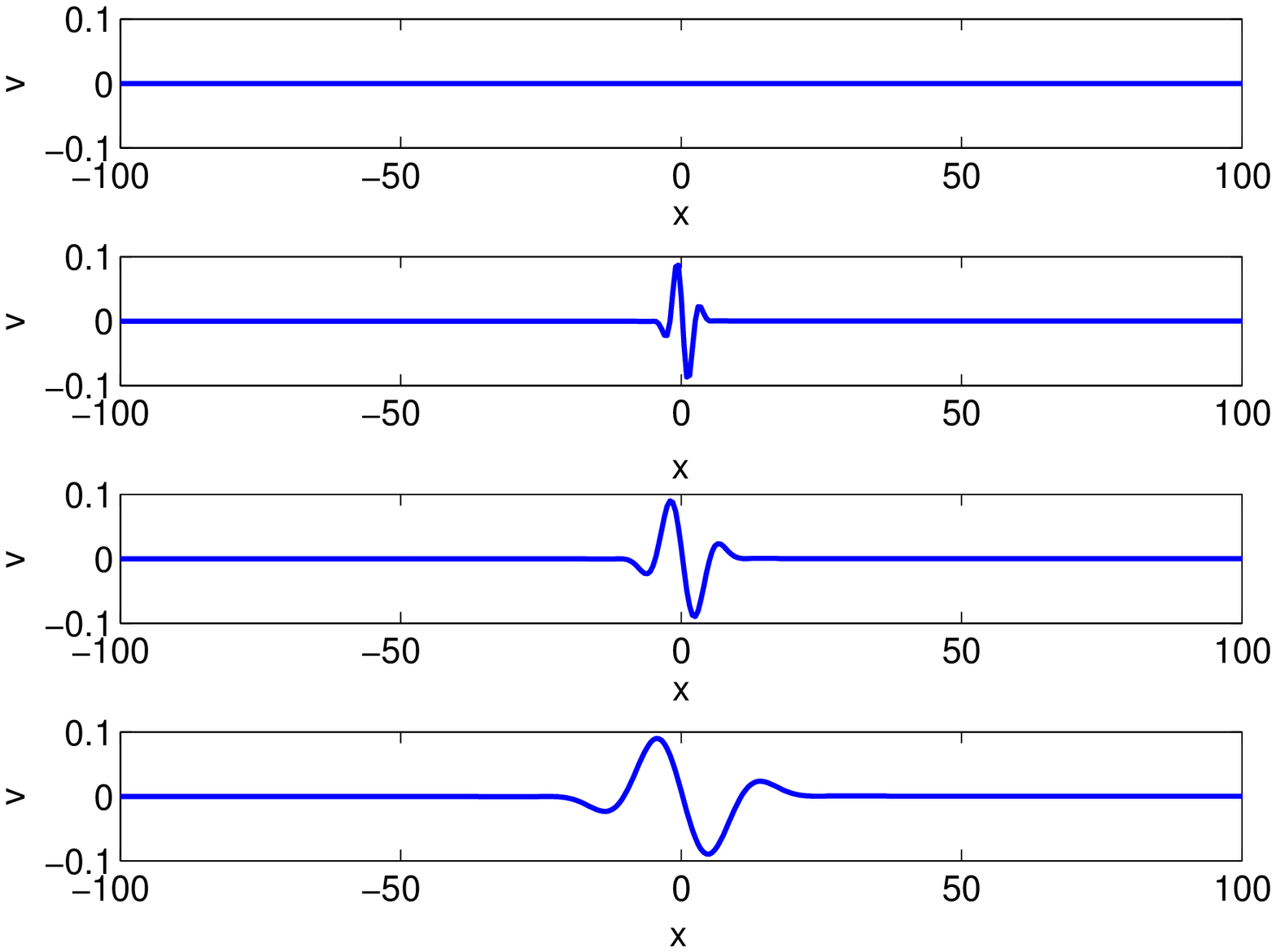}}
\caption{Cross-diffusion with $p=3$ and for $d_{11}=1,d_{22}=1.1, d_{12}=0.1,d_{21}=1$. Profiles of  (a) $u$ and (b) $v$ at times $t=0,0.25,2.5,25$.}
\label{fig_I_4}
\end{figure}

\begin{figure}[htbp]
\centering
\subfigure[]
{\includegraphics[width=0.5\textwidth]{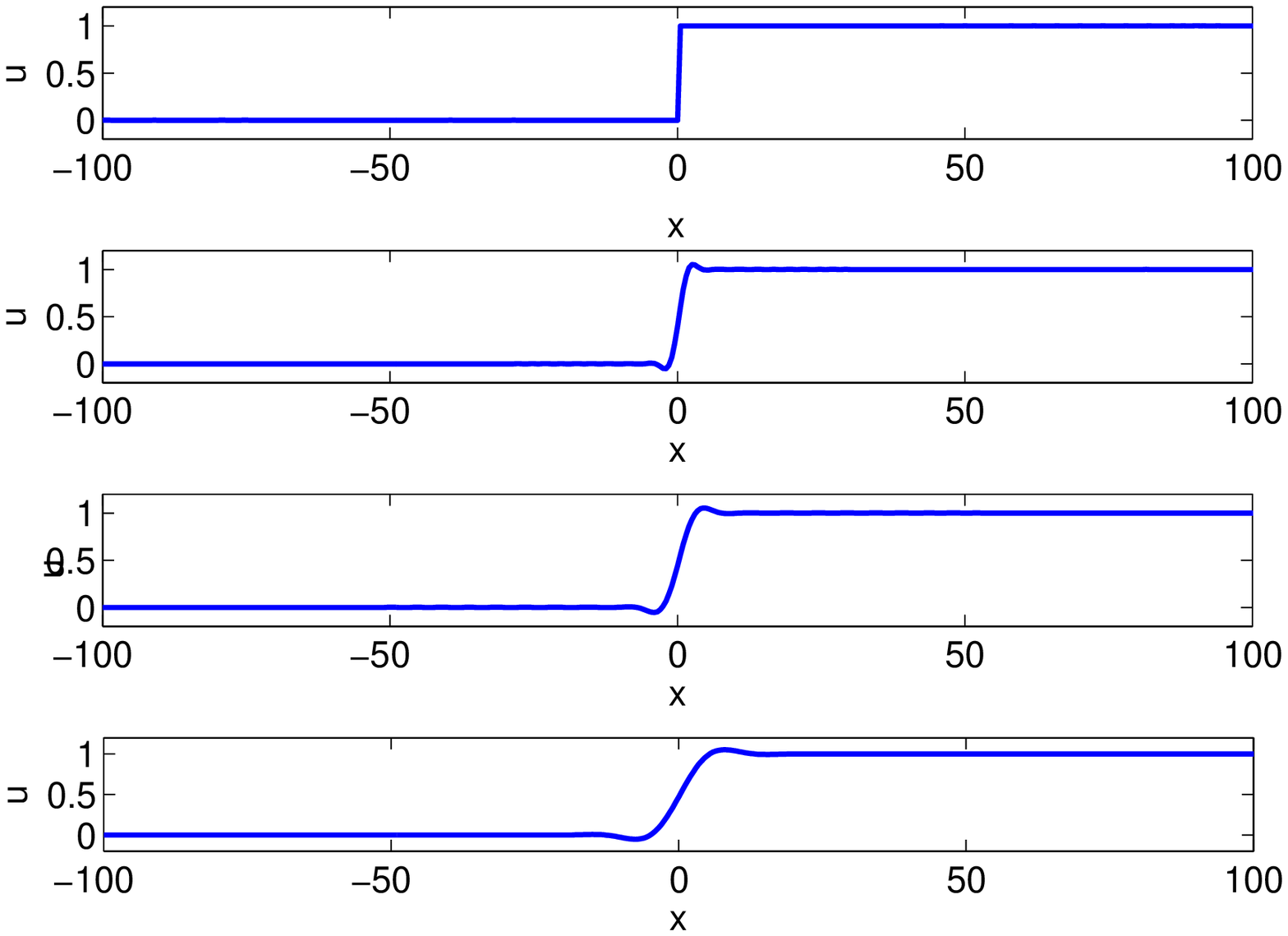}}
\subfigure[]
{\includegraphics[width=0.5\textwidth]{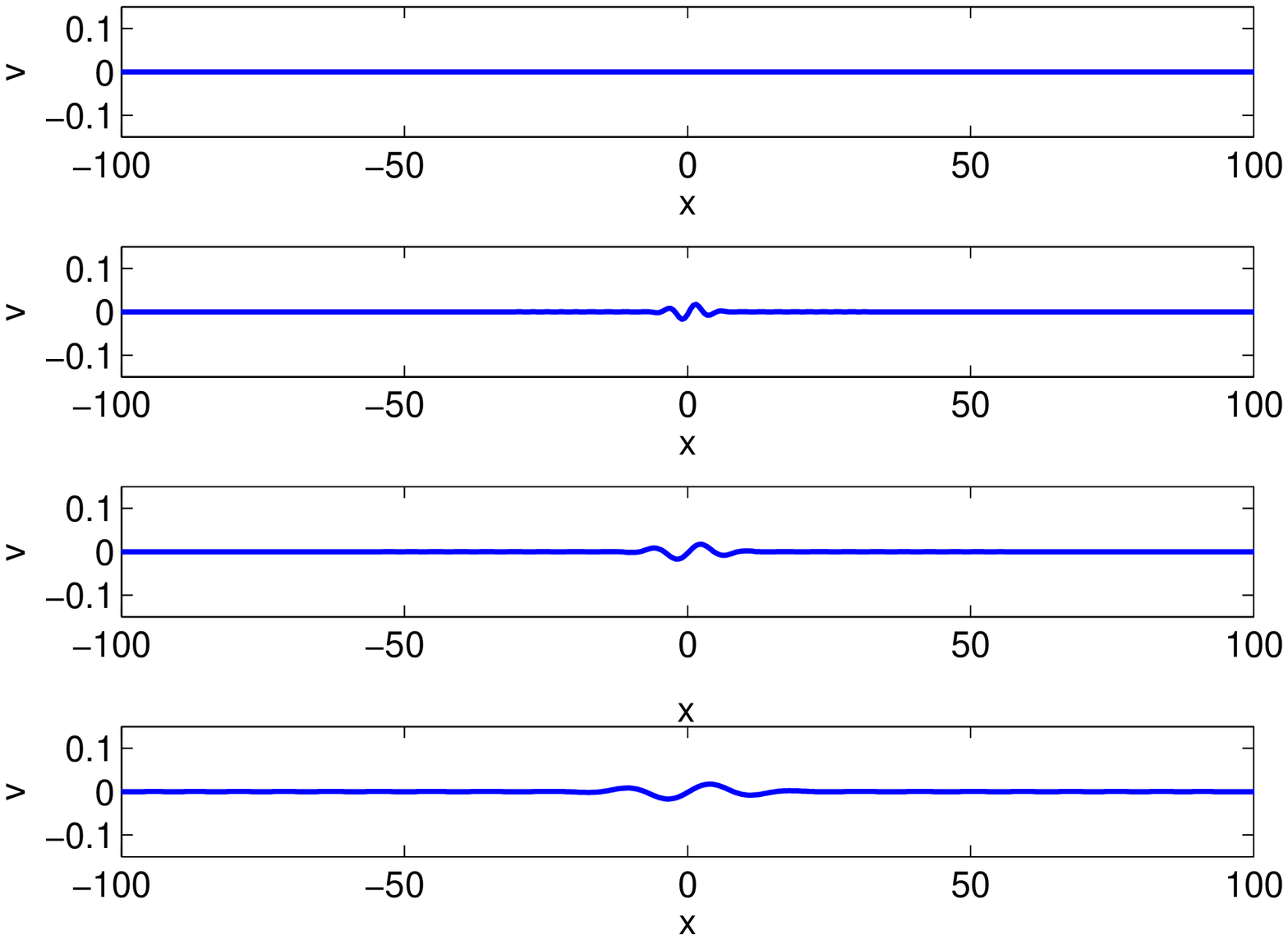}}
\caption{Cross-diffusion with $p=4$ and for $d_{11}=1,d_{22}=1.1, d_{12}=0.1,d_{21}=-0.25$. Profiles of  (a) $u$ and (b) $v$ at times $t=0,0.25,2.5,25$.}
\label{fig_I_5}
\end{figure}

\begin{figure}[htbp]
\centering
\subfigure[]
{\includegraphics[width=0.5\textwidth]{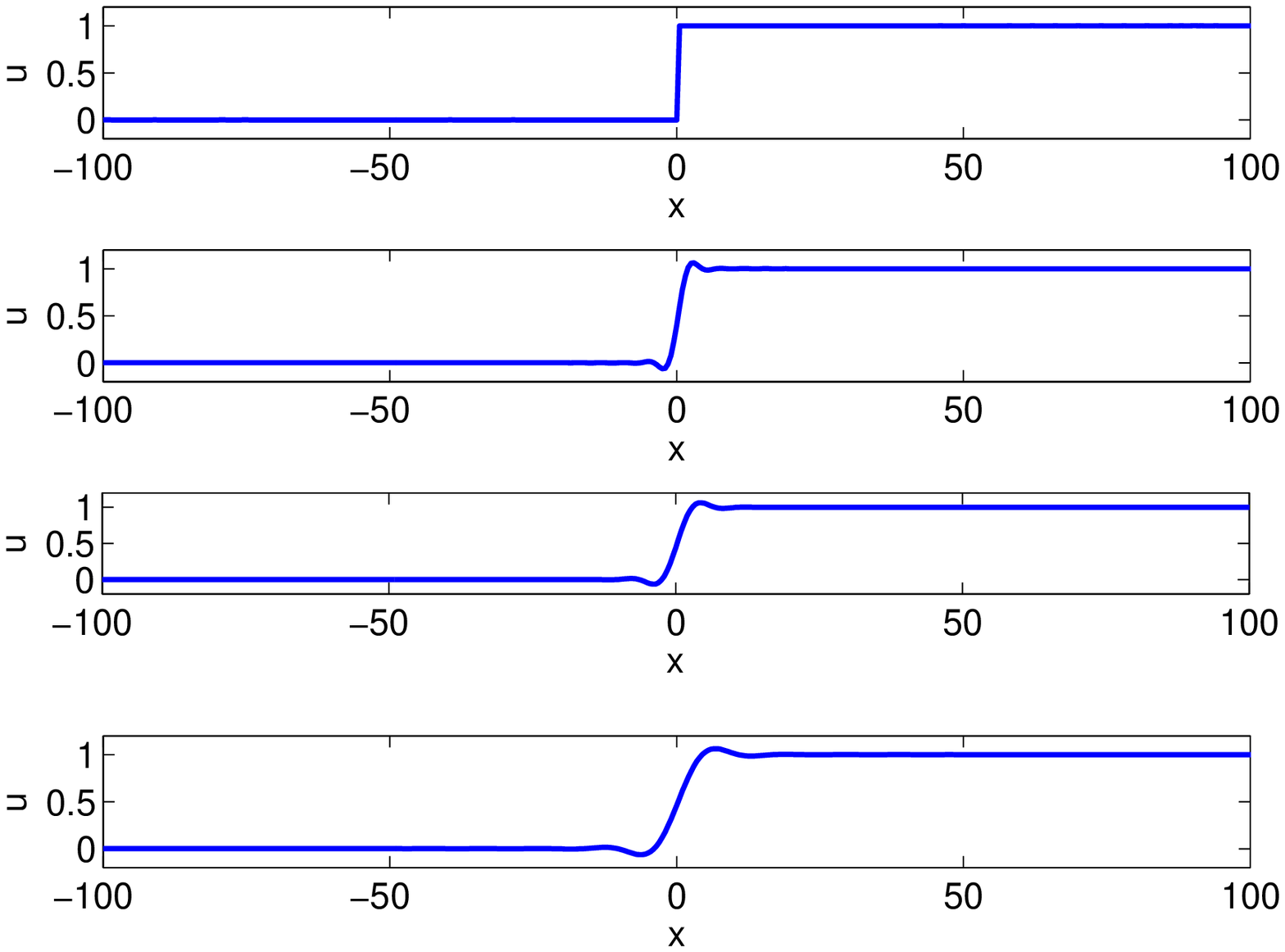}}
\subfigure[]
{\includegraphics[width=0.5\textwidth]{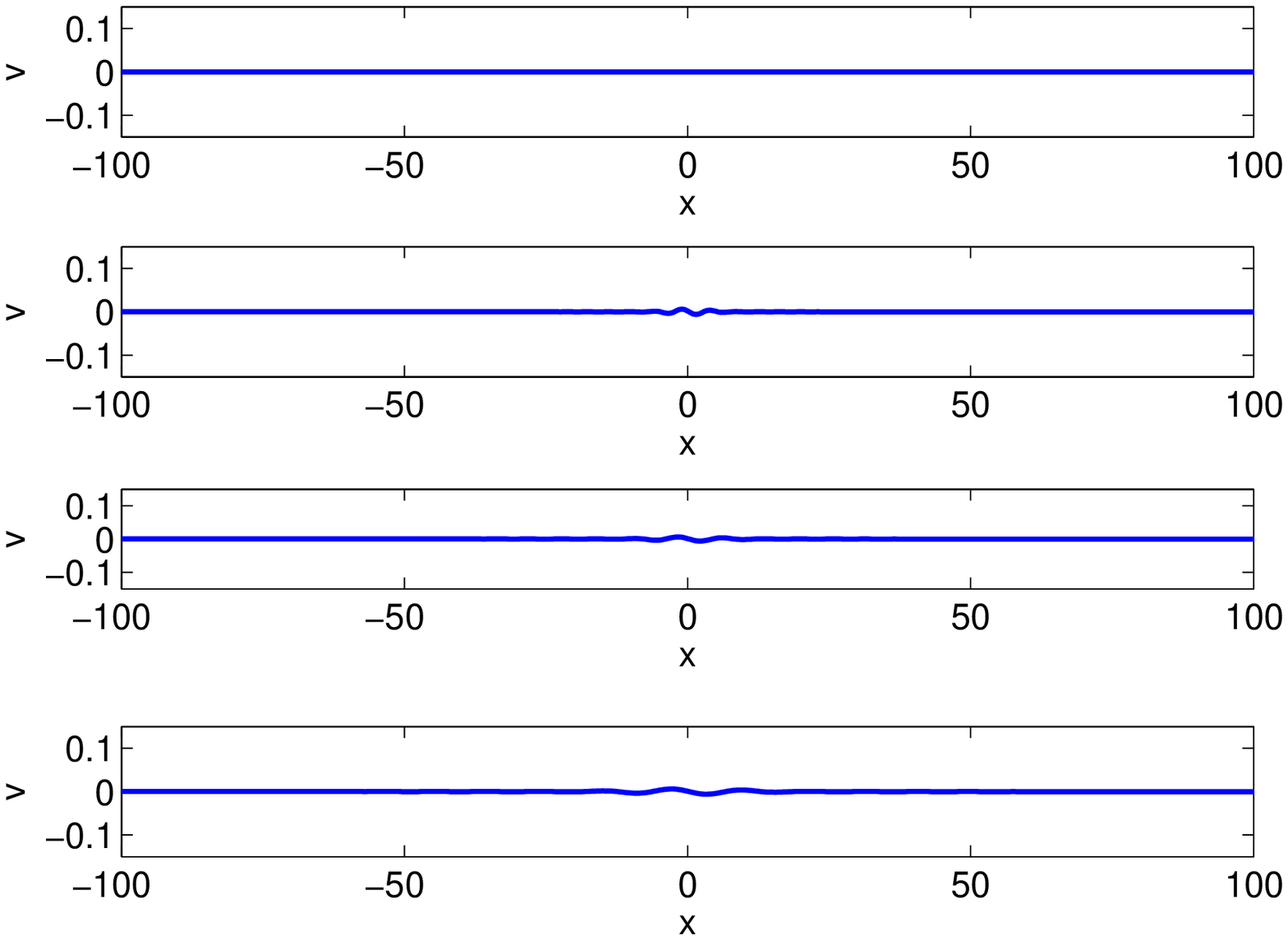}}
\caption{Cross-diffusion with $p=5$ and for $d_{11}=d_{22}=1, d_{12}=-0.1,d_{21}=0.1$. Profiles of  (a) $u$ and (b) $v$ at times $t=0,0.25,2.5,25$.}
\label{fig_I_6}
\end{figure}

The influence of the values of $p$ is first illustrated in Figures \ref{fig_I_4}-\ref{fig_I_6} for 1D signals and similar matrices to those of cases (i), (iii) and (iv) in Lemma \ref{lem_app}. Note that as $p$ grows the first component develops small oscillations at the points with less regularity.  On the other hand, the second component increases somehow the number of pulses. 
\subsubsection{Experiments in 2D}
\label{sec322}
\begin{figure}[htbp]
\centering
\subfigure[]
{\includegraphics[width=0.5\textwidth]{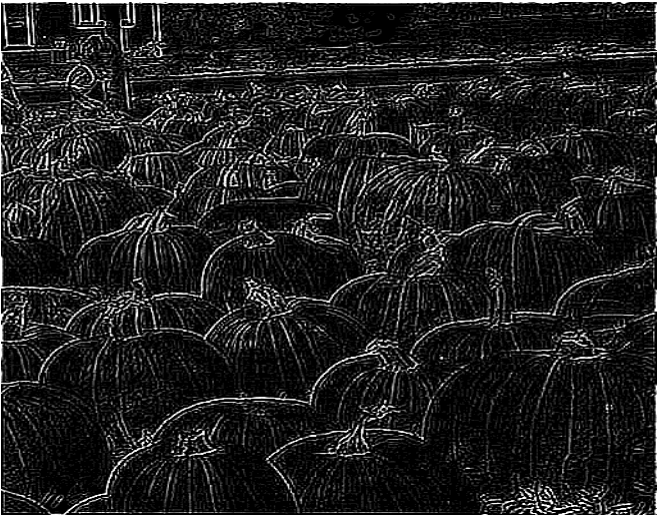}}
\subfigure[]
{\includegraphics[width=0.5\textwidth]{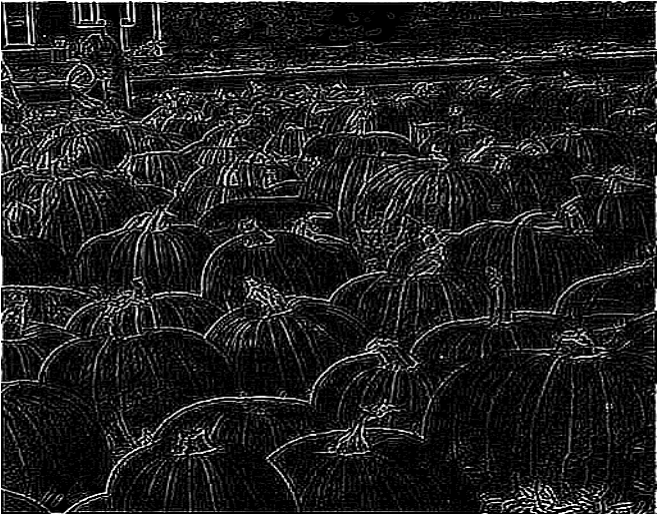}}
\caption{Generalized small theta approximation.  Second component of the solution of (\ref{cdl3}) at $t=0.1$ with $p=3, {\bf u}_{0}=(f,0)^{T}$ and (b) $d=d_{1}$, (c) $d=d_{2}$, see (\ref{cdl40}).}
\label{figR_F}
\end{figure}

As in Section \ref{sec312}, a first point to study here is the effect of $p$ on the generalized small theta approximation property. A similar experiment with $d_{1}$ and $d_{2}$ given by (\ref{cdl40}) illustrates the assignment of the second component as edge detector, under the conditions of Theorem \ref{the3} but when the operator $D$ in (\ref{cdl25}) is nonlocal. The results are shown in Figure \ref{figR_F}, corresponding to $p=3$. Compared to Figure \ref{figR_A} (for which $p=2$) no relevant differences are observed.

\begin{figure}[htbp]
\centering
\subfigure[]
{\includegraphics[width=0.5\textwidth]{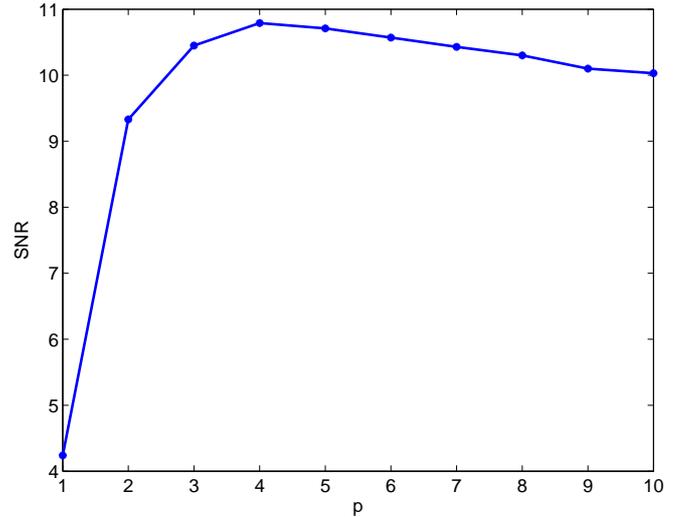}}
\subfigure[]
{\includegraphics[width=0.5\textwidth]{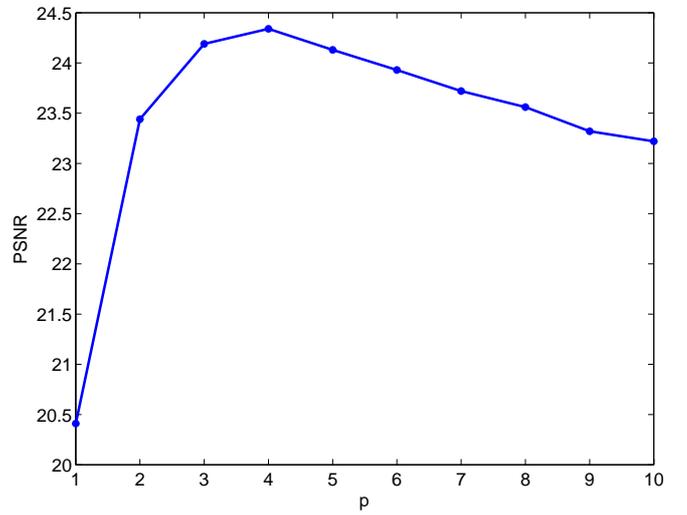}}
\caption{$SNR$ (a) and $PSNR$ (b)  vs. $p$ for a filter (\ref{cdl3}) at $t=10$ with ${\bf u}_{0}^{(1)}=(f,0)^{T}$ where $f$ is the initial noisy image affected by additive Gaussian white noise with $\sigma=40$ and $d=d_{1}$, see (\ref{cdl41}).}
\label{figR_G}
\end{figure}

\begin{figure}[htbp]
\centering
\subfigure[]
{\includegraphics[width=0.5\textwidth]{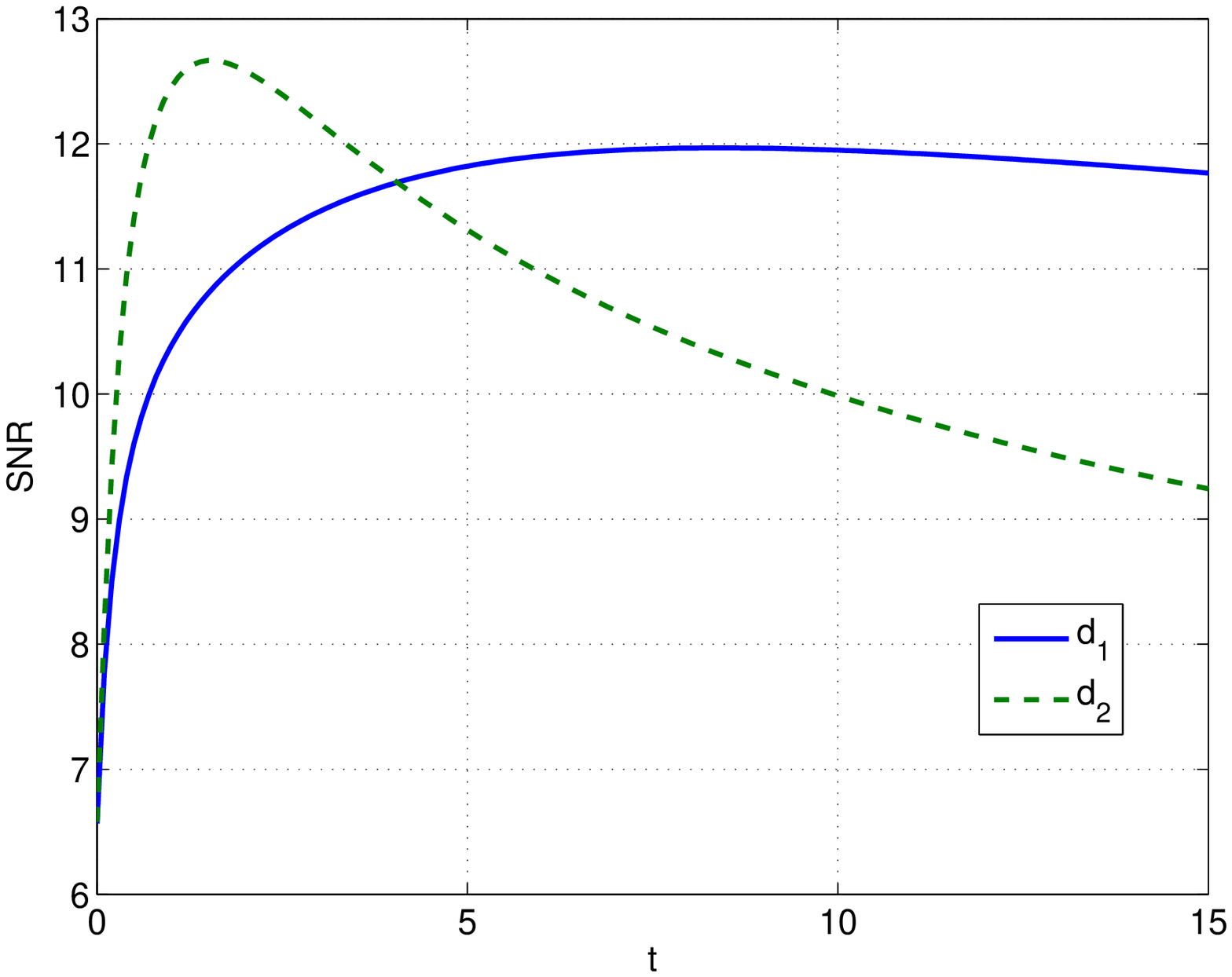}}
\subfigure[]
{\includegraphics[width=0.5\textwidth]{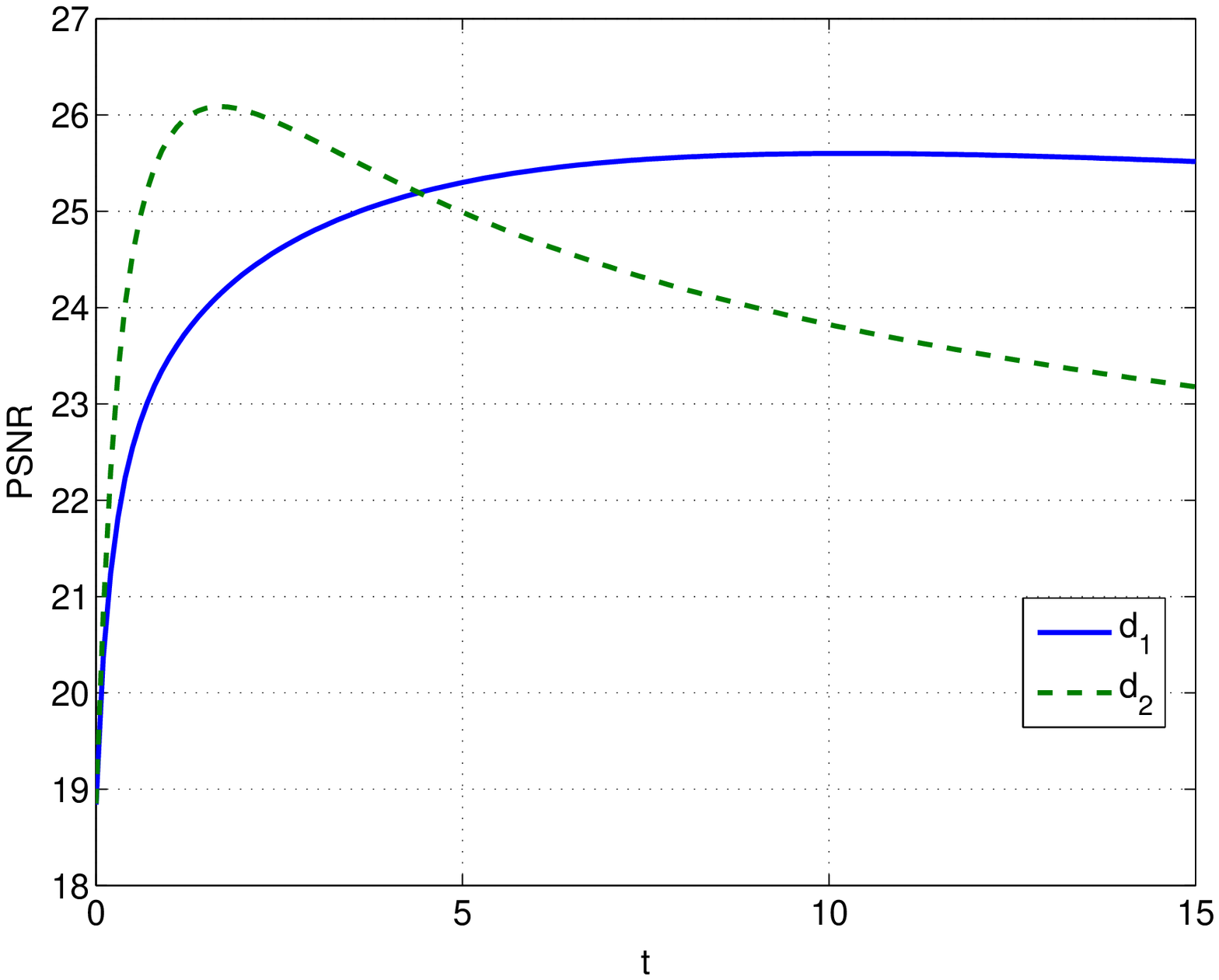}}
\caption{(a) $SNR$ vs. $t$ and (b) $PSNR$ vs. $t$ for a filter (\ref{cdl3}) with $p=3, {\bf u}_{0}=(f,0)^{T}$, $f$ in Figure \ref{figR_C}(a), {for both  $d=d_{1}$ and  $d=d_{2}$}, see (\ref{cdl41}).}
\label{figR_H}
\end{figure}

The influence of the local or nonlocal character of (\ref{cdl25}) on the quality of filtering was studied by numerical means and some of the experiments performed are shown here. The first one, displayed in Figure \ref{figR_G}, compares  the SNR and PSNR values as functions of $p$ obtained by computing (\ref{cdl3}) from a noisy image $f$ with ${\bf u}_{0}=(f,0)^{T}$,  $d=d_{1}$ given by (\ref{cdl41}) and at time $t=10$. Note that both parameters attain a maximum value by $p=4$ (for which the operator $D$ in (\ref{cdl18}) is local) while  {among the nonlocal generators, those around the values $p=3$ (in the case of PSNR) and $p=5$ (in the case of SNR) show the best results}. In all the related experiments performed, the same behaviour was observed.
As for the evolution of the filtering, a similar experiment to that of Figure \ref{figR_D} but for different values of $p$ is illustrated in Figure \ref{figR_H}. The behaviour of the SNR and PSNR values suggest that the advantages of using (\ref{cdl3}) with matrices of the type of $d_{1}$ are independent of of $p$.

\subsection{Choice of the initial distribution}
\label{sec33}
A final question in this numerical study is the influence of the initial distribution ${\bf u}_{0}$ in the small theta approximation and in filtering problems in 2D with (\ref{cdl3}). As far as the first one is concerned, note that Theorem \ref{the3} is applied for ${\bf u}_{0}=(f,0)^{T}$ or ${\bf u}_{0}=(0,f)^{T}$ where $f$ is the original image. As $s\rightarrow 0$, in all the cases
(\ref{cdl37})-(\ref{cdl39}) , $u({\bf x},t), v({\bf x},t)$ behave as
\begin{eqnarray}
u({\bf x},t)&\approx&e^{\frac{q}{2} tA}\left((1-\frac{rt}{2}A)u_{0}({\bf x})+d_{12}tAv_{0}({\bf x})\right),\nonumber\\
v({\bf x},t)&\approx&e^{\frac{q}{2} tA}\left((1+\frac{rt}{2}A)v_{0}({\bf x})+d_{21}tAu_{0}({\bf x})\right).\label{cdl42b}
\end{eqnarray}
Then the approximation (\ref{cdl42b}) (which is actually exact in the case of (\ref{cdl38})) suggests to explore, at least numerically, some other choices for ${\bf u}_{0}$. Among the ones used in our numerical experiments, by way of illustration two are considered here, namely
\begin{eqnarray}
{\bf u}_{0}^{(1)}=(f,|\nabla f|)^{T}, {\bf u}_{0}^{(2)}=(f,-|\nabla f|\Delta f)^{T}, \label{cdl43}
\end{eqnarray}
and we denote ${\bf u}_{0}^{(0)}=(f,0)^{T}$.

\begin{figure}[htbp]
\centering
\subfigure[]
{\includegraphics[width=0.5\textwidth]{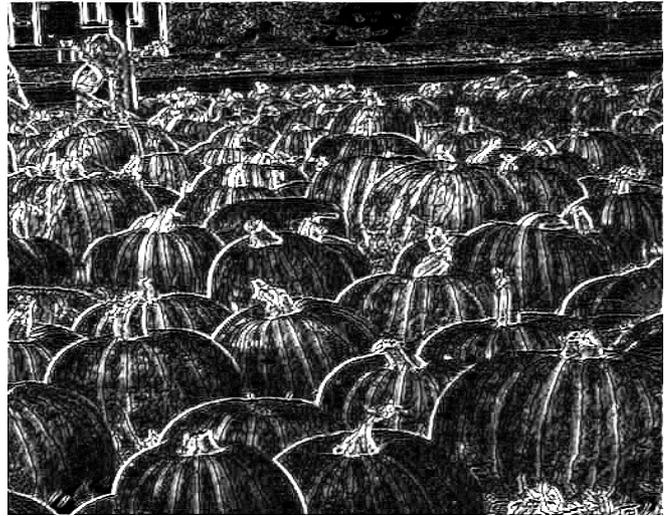}}
\subfigure[]
{\includegraphics[width=0.5\textwidth]{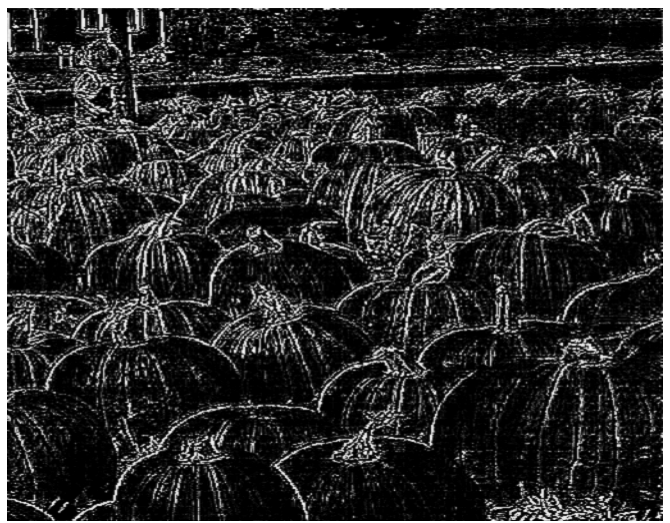}}
\caption{Generalized small theta approximation.  Second component of the solution of (\ref{cdl3}) at $t=0.1$ with $p=2,d=d_{1}$ in (\ref{cdl40}) and: (a) ${\bf u}_{0}={\bf u}_{0}^{(1)}$; (b) ${\bf u}_{0}={\bf u}_{0}^{(2)}$, see (\ref{cdl43}).}
\label{figR_I}
\end{figure}

\begin{figure}[htbp]
\centering
\subfigure[]
{\includegraphics[width=0.5\textwidth]{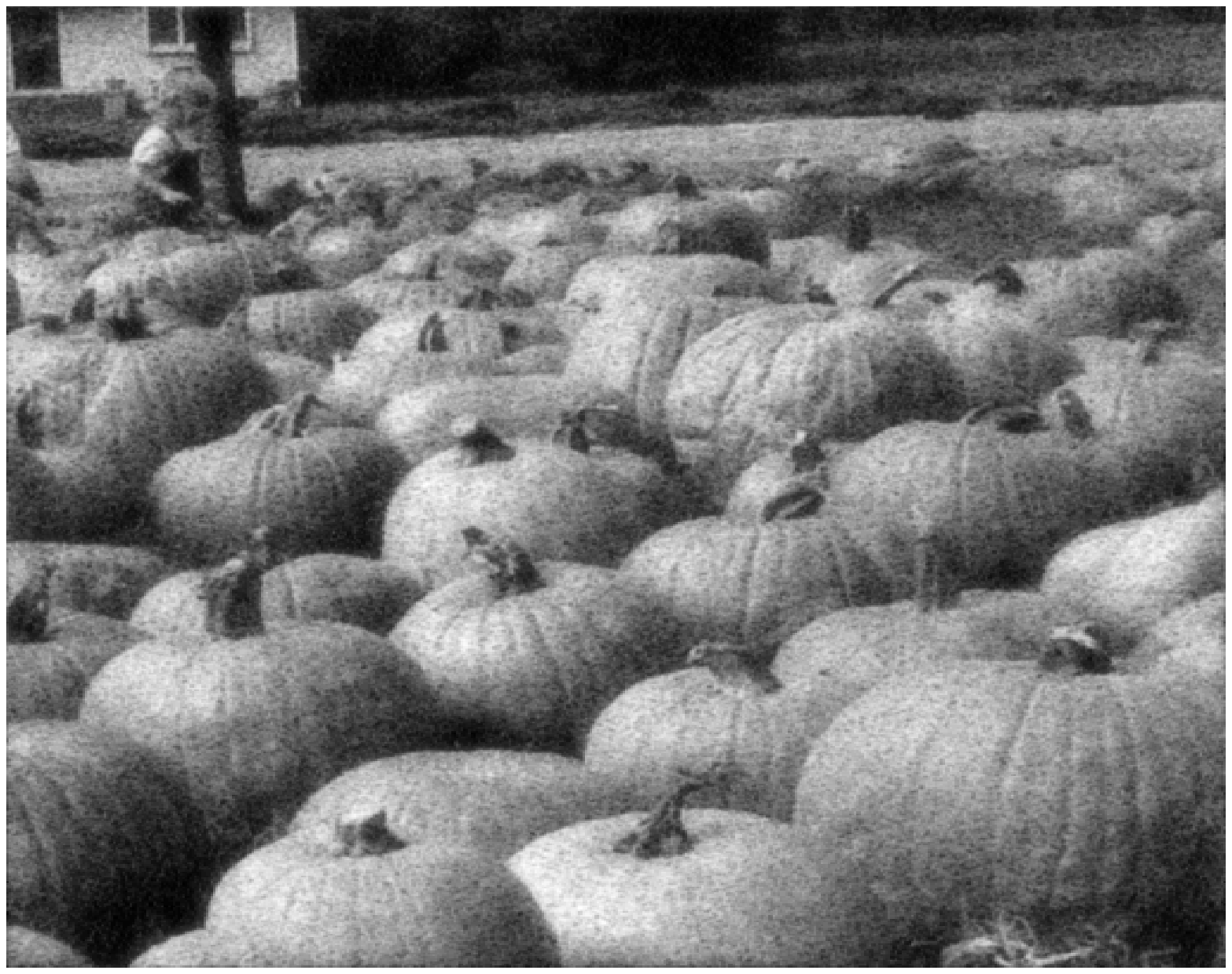}}
\subfigure[]
{\includegraphics[width=0.5\textwidth]{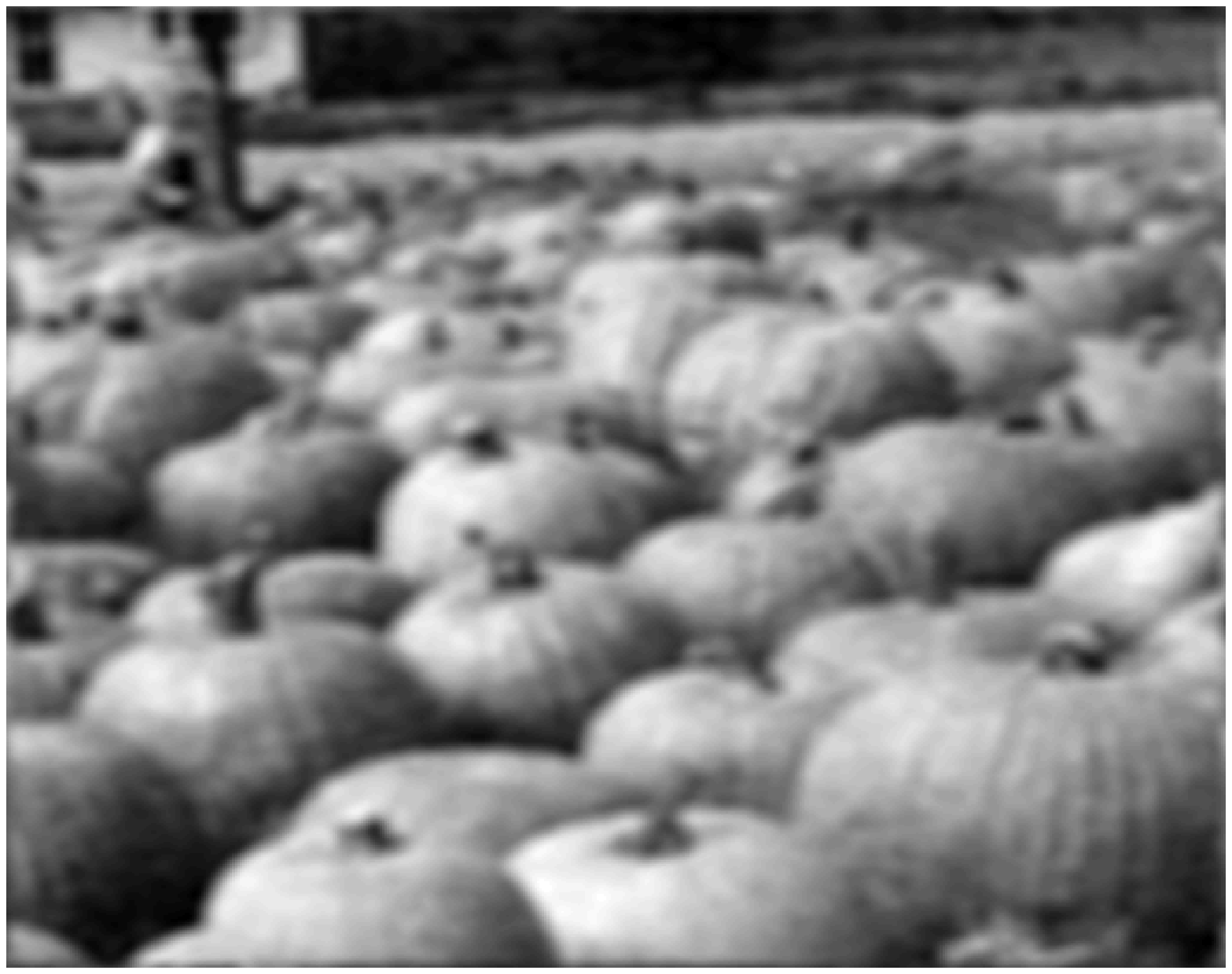}}
\caption{Image filtering problem. First component of the solution of (\ref{cdl3}) at $t=15$ with $p=2, {\bf u}_{0}={\bf u}_{0}^{(1)}$ and: (a) $d=d_{1}$, (b) $d=d_{2}$, see (\ref{cdl41}).}
\label{figR_J}
\end{figure}

\begin{figure}[htbp]
\centering
\subfigure[]
{\includegraphics[width=0.5\textwidth]{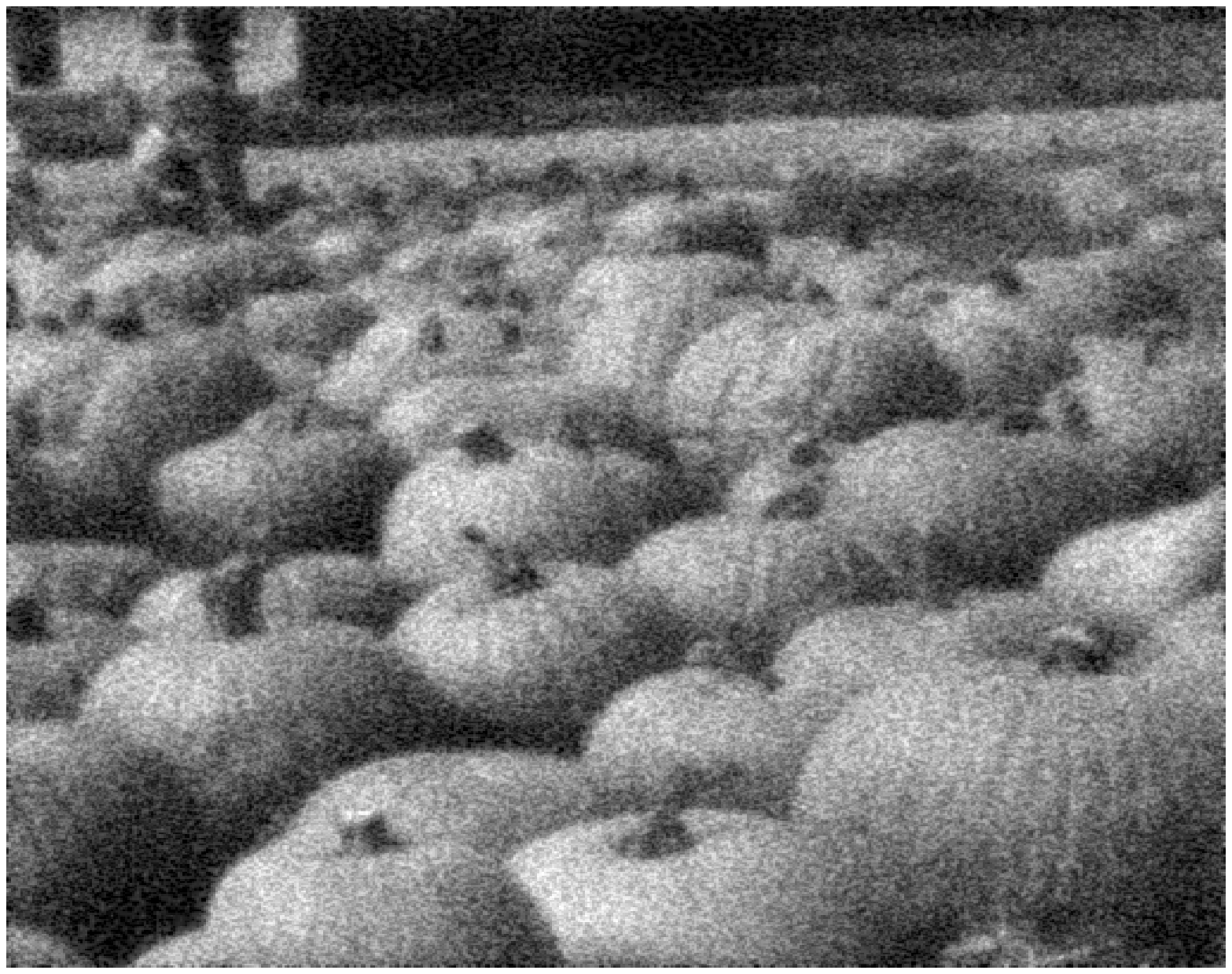}}
\subfigure[]
{\includegraphics[width=0.5\textwidth]{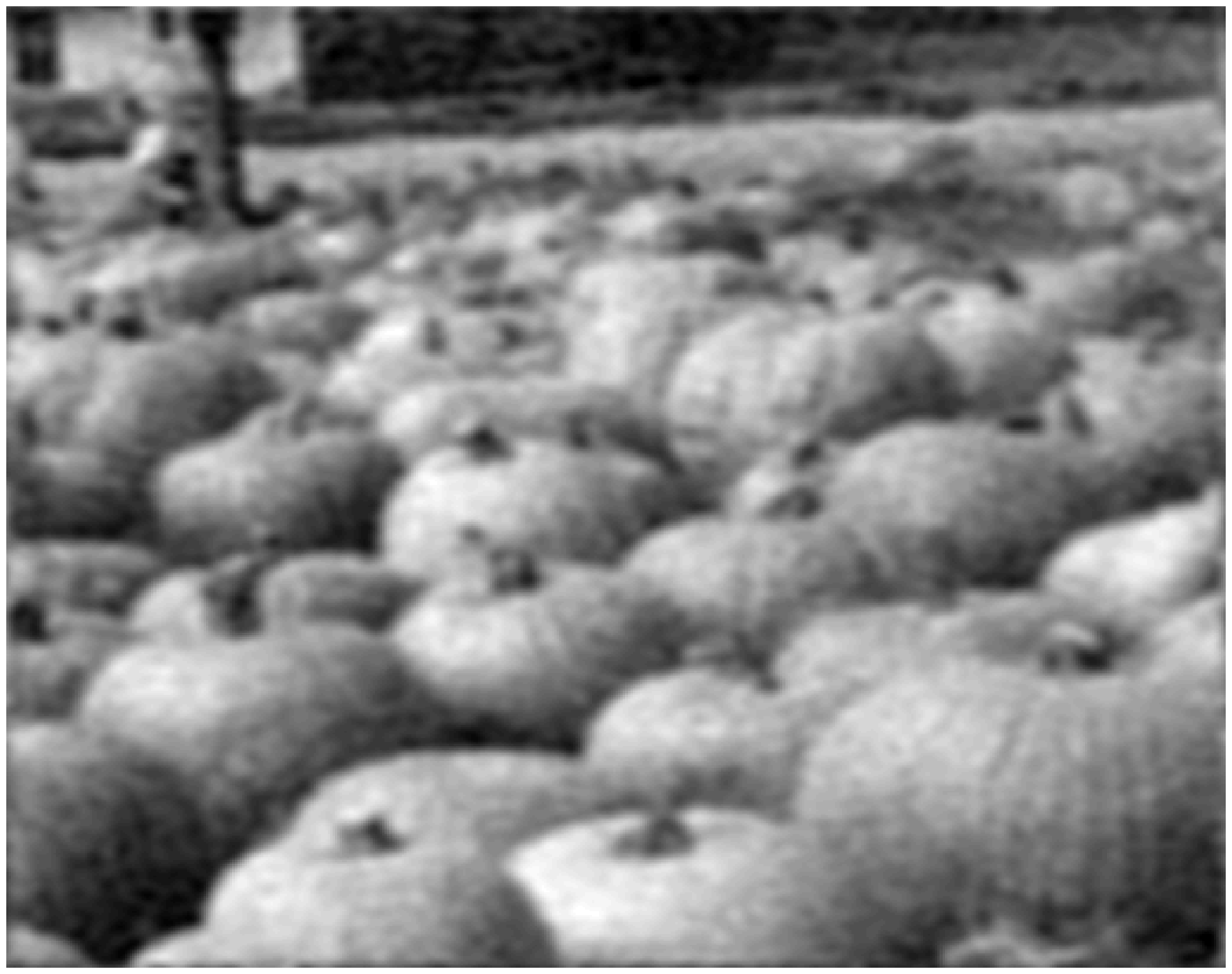}}
\caption{Image filtering problem. First component of the solution of (\ref{cdl3}) at $t=15$ with $p=2, {\bf u}_{0}={\bf u}_{0}^{(2)}$ and: (a) $d=d_{1}$, (b) $d=d_{2}$, see (\ref{cdl41}).}
\label{figR_K}
\end{figure}

By making similar experiments to that of Figure \ref{figR_A} in Section \ref{sec311}, the results are illustrated in Figure \ref{figR_I}, which corresponds to initial data given by (\ref{cdl43}).
As for the quality of filtering, the experiments presented here can be compared to those of Section \ref{sec312} for ${\bf u}_{0}={\bf u}_{0}^{(0)}$. Specifically, Figures \ref{figR_J} and \ref{figR_K} correspond to the same experiment as in Figure \ref{figR_C} but, respectively, with  ${\bf u}_{0}={\bf u}_{0}^{(1)}$ and ${\bf u}_{0}={\bf u}_{0}^{(2)}$.
In the three experiments we observe that the results do not improve those obtained  by ${\bf u}_{0}^{(0)}$ in Section \ref{sec312}. This is confirmed by a last experiment, where the original image from Figure \ref{figR_A}(a) is blurred by adding Gaussian white noise with several values of standard deviation $\sigma$. Then (\ref{cdl3}) with $p=2, d=d_{1}$ in (\ref{cdl41}) and the three different initial distributions ${\bf u}_{0}^{(j)}, j=0,1,2$ is applied. Then the SNR and PSNR values (\ref{snr}) and (\ref{psnr}) at time  $t=5$ are computed. The results are displayed in Table \ref{tavR_B} and show that the best values of the metrics are given by ${\bf u}_{0}^{(0)}$.
\begin{table}
\caption{$SNR$ (top) and $PSNR$ (bottom) values for (\ref{cdl3}) at $t=5$ with $p=2, d=d_{1}$ in (\ref{cdl41}) and ${\bf u}_{0}={\bf u}_{0}^{(j)}, j=0,1,2$.}\label{tavR_B}
\begin{center}
\begin{tabular}{cccc}
\hline
 $\sigma$&$15$&$25$&$35$\\
\hline
${\bf u}_{0}^{(0)}$&$13.31$&$11.66$&$10.05$\\
${\bf u}_{0}^{(1)}$&$12.92$&$11.24$&$9.62$\\
${\bf u}_{0}^{(2)}$&$5.19$&$4.32$&$3.74$\\
\hline
\end{tabular}
\end{center}
\begin{center}
\begin{tabular}{cccc}
 $\sigma$&$15$&$25$&$35$\\
\hline
${\bf u}_{0}^{(0)}$&$27.05$&$25.40$&$23.78$\\
${\bf u}_{0}^{(1)}$&$26.66$&$24.97$&$23.34$\\
${\bf u}_{0}^{(2)}$&$17.10$&$15.95$&$15.22$\\
\hline
\end{tabular}
\end{center}
\end{table}
\section{Conclusions and perspectives}
\label{sec4}
In the present paper linear cross-diffusion systems for image processing are analyzed. Viewed as convolution processes, those kernels satisfying fundamental scale-space properties are characterized in terms of a positive definite matrix to control the cross-diffusion and a positive parameter that determines the local character of the infinitesimal generator.
The axiomatic is based on scale invariance and 
generalizes that of the scalar case. The cross-diffusion approach, viewed as a generalization of linear complex diffusion, is shown to satisfy more general versions of the small theta approximation property to assign a role of edge detector to one of the components.

In a second part, a numerical study of comparison with kernels is made. This can be considered as a first approach, by computational means, to the performance of linear cross-diffusion models, to be analyzed in a more exhaustive way in future works.
The numerical experiments, performed for one- and two-dimensional signals, show the influence of the choice of the initial distribution of the image in a vector of two components, as well as of the matrix of the kernel on the behaviour of the {filtering} process by cross-diffusion. The numerical results suggest that suitable choices of the positive definite matrix give a delay of blurring which can also be useful to a better identification of the edges and that is independent of the local or nonlocal character of the infinitesimal generator. Additionally, other values of the initial distribution, different from the ones for which the generalized small theta approximation holds, do not improve the results {in our experiments}.

The present paper will be continued in a natural way by the introduction of nonlinear cross-diffusion models and the study of their behaviour in image restoration (see Ara\'ujo et al. \cite{ABCD2016_II}).
\begin{acknowledgements}
This work was supported by  Spanish Ministerio de Econom\'{\i}a y Competitividad under the Research Grant MTM2014-54710-P.
A. Ara\'ujo and S. Barbeiro were also supported by the Centre for Mathematics of the University of Coimbra -- UID/MAT/00324/2013, funded by the Portuguese Government through FCT/MCTES and co-funded by the European Regional Development Fund through the Partnership Agreement PT2020.
\end{acknowledgements}

\section*{Appendix}

Note that the matrix $d$ in (\ref{cdl25}) must be positive definite, but not necessarily symmetric. (This means that ${\bf x}^{T}d{\bf x}>0$ for all ${\bf x}\neq 0$ or, equivalently, if the symmetric part $(d+d^{T})/2$ is positive definite.) This implies that the real part of each eigenvalue is positive. In terms of the entries of $d$, the positive definite character requires two conditions
\begin{eqnarray}
d_{11}>0,\quad 4d_{11}d_{22}-(d_{12}+d_{21})^{2}>0,\label{cdl26}
\end{eqnarray}
or equivalently, being
\begin{eqnarray}
\lambda_{\pm}&=&\frac{1}{2}\left(d_{11}+d_{22}\pm\sqrt{s}\right), \nonumber\\
s&=&r^{2}+4d_{12}d_{21},\quad r=d_{22}-d_{11}, \label{cdl27}
\end{eqnarray}
the eigenvalues of $d$, then $Re(\lambda_{\pm})>0$.
The following spectral analysis of a positive definite matrix $d$ is used in several results of the paper. The proof is straightforward by using the standard Jordan reduction theory.
\begin{lemma}
\label{lem_app}
Assume that $d=(d_{ij})_{i,j=1,2}$ is a $2\times 2$ positive definite matrix with eigenvalues $\lambda_{\pm}$ given by (\ref{cdl27}). Then one of the following cases holds: 
\begin{itemize}
\item[(i)] If $s>0$ then $\lambda_{+}>\lambda_{-}>0$ and $d=P\Lambda P^{-1}$ with
\begin{eqnarray*}
\Lambda=\begin{pmatrix}\lambda_{+}&0\\0&\lambda_{-}\end{pmatrix},
\end{eqnarray*}
and $P$ depends on the nondiagonal entries of $d$. Explicitly:
\begin{eqnarray*}
P=\begin{pmatrix}d_{12}&\frac{r+\sqrt{s}}{2}\\d_{12}&\frac{r-\sqrt{s}}{2}\end{pmatrix}, \quad P=\begin{pmatrix}\frac{-r+\sqrt{s}}{2}&d_{21}\\\frac{-r-\sqrt{s}}{2}&d_{21}\end{pmatrix},\quad P=I,
\end{eqnarray*}
when, respectively, $d_{12}\neq 0$, $d_{21}\neq 0$ and $d_{12}=d_{21}=0$, {where $I$ is the $2\times 2$ identity matrix}.
\item[(ii)] If $s=0$ and $d$ is diagonalizable, then $\alpha=\lambda_{+}=\lambda_{-}>0$ and $d=P\Lambda P^{-1}$ with $\Lambda=\alpha I, P=I$.
\item[(iii)] If $s=0$ and $d$ is not diagonalizable, then $\alpha=\lambda_{+}=\lambda_{-}>0$ and $d=P\Lambda P^{-1}$ with \begin{eqnarray*}
\Lambda=\begin{pmatrix}\alpha&1\\0&\alpha\end{pmatrix},
\end{eqnarray*}
and 
\begin{eqnarray*}
P=\begin{pmatrix}d_{12}&\frac{r}{2}\\0&1\end{pmatrix}, \quad \begin{pmatrix}-\frac{r}{2}&d_{21}\\1&0\end{pmatrix},
\end{eqnarray*}
when, respectively, $d_{12}\neq 0$ and $d_{21}\neq 0$.
\item[(iv)] If $s<0$ then $\lambda_{\pm}=\nu\pm i\mu$ with $\nu=(d_{11}+d_{22})/2>0, \mu=\sqrt{-s}/2\neq 0$ and $d=P\Lambda P^{-1}$ with
\begin{eqnarray*}
\Lambda=\begin{pmatrix}\nu&-\mu\\\mu&\nu\end{pmatrix},
\end{eqnarray*}
and
\begin{eqnarray*}
P=\begin{pmatrix}d_{12}&\frac{r}{2}\\0&-\mu\end{pmatrix}, \quad \begin{pmatrix}-\frac{r}{2}&d_{21}\\-\mu&0\end{pmatrix},
\end{eqnarray*}
when, respectively, $d_{12}\neq 0$ and $d_{21}\neq 0$.
\end{itemize}

\end{lemma}
Note that in the case (iv) $P$ is formed by using the invariant subspace generated by the real and imaginary parts of a basis of eigenvectors of $d$.

\end{document}